\documentclass[11pt,reqno]{amsart}
\usepackage[left=3cm,top=2.5cm,right=3cm,bottom=2.5cm]{geometry}
\usepackage{amsmath,amsfonts}
\usepackage{amssymb,amsthm}
\usepackage{graphicx}
\usepackage[latin1]{inputenc}
\usepackage[T1]{fontenc}
\usepackage[english]{babel}

\theoremstyle{plain}

\newtheorem{Teo}{Theorem}[section]
\newtheorem{Lema}[Teo]{Lemma}
\newtheorem{Prop}[Teo]{Proposition}

\theoremstyle{remark}
\newtheorem{Example}[Teo]{Example}

\def\w{w}
\def\r{\tilde{r}}
\def\ww{w^{\ast}}
\newcommand{\etal}{\textsl{et al.}}

\begin{document}

\pagestyle{myheadings} 
\title{Characterizing trees with large Laplacian energy}

\author[E. Fritscher]{Eliseu Fritscher}\email{\tt eliseu.fritscher@ufrgs.br}

\author[C. Hoppen]{Carlos Hoppen}\email{\tt choppen@ufrgs.br}

\author[I. Rocha]{Israel Rocha} \email{\tt israel.rocha@ufrgs.br}

\author[V. Trevisan]{Vilmar Trevisan} \email{\tt trevisan@mat.ufrgs.br}

\address{Instituto de Matem\'atica, UFRGS -- Avenida Bento Gon\c{c}alves, 9500, 91501--970 Porto Alegre, RS, Brazil}
\thanks{C. Hoppen acknowledges the support of FAPERGS (Proc.~11/1436-1). V. Trevisan was partially supported by CNPq (Proc. 309531/2009-8 and 473815/2010-9) and FAPERGS (Proc. 11/1619-2)}

\begin{abstract}
We investigate the problem of ordering trees according to their Laplacian energy. More precisely, given a positive integer $n$, we find a class of cardinality approximately $\sqrt{n}$ whose elements are the $n$-vertex trees with largest Laplacian energy. The main tool for establishing this result is a new upper bound on the sum $S_k(T)$ of the $k$ largest Laplacian eigenvalues of an $n$-vertex tree $T$ with diameter at least four, where $k \in \{1,\ldots,n\}$.
\end{abstract}

\maketitle

\section{Introduction}

In this paper, we address a particular instance in the wide class of extremal problems in spectral graph theory. Problems of this type include, for instance, finding the extremal value of some spectral parameter over a class of graphs, characterizing  the elements of this class that achieve this extremal value, or ordering the elements in this class according to the value of this parameter.

A great deal of research has been done in this general framework. For results until $1988$, we refer the reader to  Cvetkovi\'{c} et al.~\cite{Cve88}. Regarding more recent work, perhaps the most studied parameter in this context is the algebraic connectivity, and the multitude of extremal results about it is illustrated by the work of Abreu~\cite{Abr07} and by the references therein. Concerning other parameters, we may cite, for example, Wang and Kang \cite{WK09} and Guo \cite{G08a}, who ordered graphs by energy, and Zhang and Chen~\cite{ZH02}, Yu, Lu and Tian~\cite{YLT05}, Lin and Guo~\cite{LG06}, Belardo, Li Marzi and Simi\`{c}~\cite{BMS07}, and Guo~\cite{SGG08}, who gave orderings based on adjacency and Laplacian spectral indices. Moreover, Rada and Uzc\'ategui~\cite{RU05} gave an ordering of chemical trees according to their Randi\`{c} index,  He and Li~\cite{HL11} ordered trees by their Laplacian coefficients, while Xu \cite{X11} determined extremal trees relative to the Harary index.

The parameter of interest in this work is the Laplacian energy of a graph, which was first introduced by  Gutman and Zhou~\cite{GZ06}. Given a graph $G$ on the vertex set $V=\{v_1,\ldots,v_n\}$, the \emph{Laplacian matrix} of $G$ is given by $L=D-A$, where $D$ is the diagonal matrix whose entry $(i,i)$ is equal to the degree of $v_i$ and $A$ is the adjacency matrix of $G$. The \emph{Laplacian spectrum} of $G$ is defined as the set of eigenvalues of $L$, which we shall denote $\mu_1 \geq \mu_2 \geq \cdots \geq \mu_n=0$. With this, the \emph{Laplacian energy} of $G$ is given by
$$LE(G)=\sum_{i=1}^{n} |\mu_i-\overline{d}|,$$
where $\overline{d}$ is the average degree of $G$.

A natural extremal question involving the Laplacian energy of a graph is to ask for the largest value of this parameter among all $n$-vertex graphs, and for the $n$-vertex graph (or for the family of such graphs) whose Laplacian energy achieves this extremal value. However, this seems to be a hard problem, which, to the best of our knowledge, is still open. More is known when this problem is restricted to some particular classes of $n$-vertex graphs. For instance, Radenkovi\'c and Gutman~\cite{RG07} ran computational experiments to investigate the correlation between the energy and the Laplacian energy of a graph and, although they have found evidence that the parameters have a similar dependence on the number of vertices and on the number of edges, the correlation breaks down when the classes consist of graphs with the same number of vertices and the same number of edges. Perhaps surprisingly, after computing the energy and the Laplacian energy of a large family of trees up to 14 vertices, Radenkovi\'c and Gutman conjectured that the $n$-vertex tree with minimum energy has maximum Laplacian energy and vice-versa. The former has been confirmed in~\cite{FHRT11}, where the current authors proved that the star $\mathcal{S}_n$ is the $n$-vertex tree with maximum Laplacian energy.

Even more is known in the case of $n$-vertex trees whose diameter is equal to three, that is, trees such that the largest distance between two of their vertices is equal to three. A tree in this class can be viewed as the union of two stars with an edge between their centers. The endpoints of this edge are called the two ends of the tree. Every tree with diameter three may be written in the form $T(a,b)$, where $n=a+b+2$ is the number of vertices and $a,b\geq 1$ determine the number of leaves incident with each end of $T(a,b)$. For fixed integers $a \geq b \geq 1$, consider
the class $\mathcal{T}(a,b)=\{T(a+k,b-k) \colon k=0, \ldots, b-1\}$. Trevisan, Carvalho, Del-Vecchio and Vinagre~\cite{TCDV10} have shown that the Laplacian energy of the elements of
$\mathcal{T}(a,b)$ decreases as $k$ gets larger.
\begin{Teo}\label{orderTrevisan}\cite{TCDV10}
The Laplacian energy of $T(a+k,b-k)$ is a strictly decreasing function of $k$, for $k=0, \ldots, b-1$.
\end{Teo}

In this paper, we focus on the problem of ordering the full class of $n$-vertex trees according to their Laplacian energy. More precisely, for any fixed $n$, we find a list of roughly $\sqrt{n}$ trees on $n$ vertices with largest Laplacian energy. With the exception of the $n$-vertex star, we prove that all trees in this class have diameter three, so that the actual order is inherited from Theorem~\ref{orderTrevisan}.  In particular, this paper settles a conjecture of Radenkovi\'{c} and Gutman~\cite{RG07} regarding the $n$-vertex trees with second, third and fourth largest Laplacian energy, which are, respectively, $T(\lceil (n-2)/2 \rceil, \lfloor (n-2)/2 \rfloor)$, $T(\lceil (n-2)/2 \rceil+1, \lfloor (n-2)/2 \rfloor-1)$ and $T(\lceil (n-2)/2 \rceil+2, \lfloor (n-2)/2 \rfloor-2)$.
\begin{Teo}\label{order}
Let
\begin{equation}
f(n)=\begin{cases} 1+\lfloor\sqrt{n-3}\rfloor, & \text{ if $n$ is even;} \\
\lfloor\frac{3}{2}+\sqrt{n-3}\rfloor, & \text{ if $n$ is odd and $n=p^2-p+3$ for some $p \in \mathbb{N}$;} \\
\lfloor\frac{1}{2}+\sqrt{n-3}\rfloor, & \text{ if $n$ is odd and cannot be written as before.}\end{cases}
\end{equation}
Among all trees with $n\geq 6$ vertices, the trees with second through $(f(n)+1)$-st largest Laplacian energy are, respectively, the trees $T(\left\lceil\frac{n-2}{2}\right\rceil+k,\left\lfloor\frac{n-2}{2}\right\rfloor-k)$, where
$0\leq k \leq f(n)-1$.
\end{Teo}
This naturally leads to the question of whether the Laplacian energy of a tree with diameter three is always larger than its counterpart for a tree with larger diameter. However, we show that this is not the case by providing a counterexample for every $n\geq 16$ (see Theorem~\ref{thm_counter}).

The main tool for proving Theorem~\ref{order} is a refinement of Theorem~1.1 in~\cite{FHRT11}, which establishes the following general upper bound on the sum of the $k$ largest Laplacian eigenvalues of an $n$-vertex tree $T$:
\begin{equation}\label{eq_bound}
S_k(T):=\sum_{i=1}^k \mu_i \leq n-2+2k-\frac{2k-2}{n}.
\end{equation}
Moreover, it follows from~\cite{FHRT11} that this upper bound is tight when $k=1$ and $T$ is a star, and that it cannot be improved by subtracting $1/n$, even if we consider trees with diameter at least three. In this work, we show that this upper bound may be improved by $2/n$ for trees with diameter at least four.
\begin{Teo}{\label{teo1}} Every tree $T$ with $n\geq6$ vertices and diameter greater
than or equal to 4 satisfies
\begin{equation}\label{new_bound} S_k(T)\leq n-2+2k-\frac{2k}{n}. \end{equation}
\end{Teo}

The upper bounds in~\eqref{eq_bound} and~\eqref{new_bound} have a more general counterpart in the form of Brouwer's Conjecture~\cite{Brobook}, which states that, given a graph $G=(V,E)$ with $n$ vertices and an integer $k \in \{1,\ldots,n\}$, the sum $S_k(G)$ satisfies
\begin{equation}\label{conj_brouwer}
S_k(G) \leq |E| +{k + 1 \choose 2}.
\end{equation}
For $k = 1$, the conjecture follows from the well-known inequality $\mu_1(G) \leq |V(G)|$, and the cases $k = n$ and $k = n - 1$ are also straightforward. Haemers, Mohammadian and Tayfeh-Rezaie~\cite{HMT10} have proved this conjecture for all graphs $G$ and $k=2$, and for all trees $T$ and arbitrary values of $k$. More recently, Du and Zhou~\cite{DZ} (see also Wang, Huang and Liu~\cite{Wang}) showed that the conjecture holds for unicyclic ($1$-cyclic) and bicyclic ($2$-cyclic) graphs, among others, where a graph is said to be \emph{$c$-cyclic} if it has $n-1+c$ edges. Using the same idea, we may use our new upper bound~\eqref{new_bound} to derive a new upper bound on $S_k(G)$ for $n$-vertex graphs with diameter at least four (see Theorem~\ref{teo_cyclic}), which has the following consequence in terms of Brouwer's Conjecture.
\begin{Teo}\label{teo_brouwer}
Let $G$ be a graph with $n\geq 6$ vertices and diameter $d \geq 4$. The inequality~\eqref{conj_brouwer} is satisfied for
$$\frac{3n-4+\sqrt{9n^2-24n+16+8e(G)n^2-8n^3}}{2n} \leq k \leq n.$$
\end{Teo}
For the class of graphs under consideration, this is a slight improvement on Corollary~3.2 in~\cite{Wang}; as it turns out, the lower bound on $k$ given in~\cite{Wang} is the same but for the linear term inside the square root, which was equal to $-8n$. However, this does provide new examples for which~\eqref{conj_brouwer} is verified. For instance,  if $n=9$ and $e=15$, the new bound ensures that~\eqref{conj_brouwer} holds for $5\leq k \leq 9$ as opposed to $6 \leq k\leq 9$. More generally, if $n=11+2j$ and $e=(6+j)(7+j)/2$, where $j \in \mathbb{N}$, the new bound ensures that~\eqref{conj_brouwer} holds for $j+5\leq k$ as opposed to $j+6 \leq k$.

The paper is organized as follows. In Section~\ref{sec_consequences} we present the main consequences of Theorem~\ref{teo1}, including the proof of Theorem~\ref{order}. The actual proof of Theorem~\ref{teo1} is the subject of the subsequent sections: an overview of the proof, as well as preliminary results and definitions, lie in Section~\ref{sec_overview}, while Section~\ref{sec_particular} contains the proof of this theorem for a special family of trees, and Section~\ref{sec_general} deals with the proof of the general case.

\section{Consequences of Theorem~\ref{teo1}}\label{sec_consequences}

In the present section, we shall see that the new bound on the sum of the largest Laplacian eigenvalues given in Theorem~\ref{teo1} implies Theorem~\ref{order}. Given an integer $n \geq 6$, this theorem provides a family of cardinality $f(n)$ consisting of trees with diameter three whose Laplacian energies rank from second largest to $(f(n)+1)$-st largest among all $n$-vertex trees. To prove this result, we shall compare the Laplacian energies of an $n$-vertex tree $T(a,b)$ with diameter three and of an $n$-vertex tree with diameter at least four. A useful step in this direction is to relate the Laplacian energy of a graph $G$ with the sum of its Laplacian eigenvalues that are larger than the average degree of $G$.

\begin{Prop}\label{prop_1} Let $G$ be a graph with $\sigma$ Laplacian eigenvalues that are larger than the average degree $\bar{d}$ of $G$. The Laplacian energy of $G$ is given by
$$LE(G)=2 S_{\sigma}(G)-2\sigma\overline{d}.$$
\end{Prop}
The proof of this result relies on straightforward algebraic manipulations, and is therefore omitted (see \cite[Proof of Theorem~1.3]{FHRT11} for a detailed account). Clearly, if we wish to use Proposition~\ref{prop_1} to compare the Laplacian energies of two graphs, we need information about the sum of their largest Laplacian eigenvalues and about the number of eigenvalues that are larger than the average. For trees with diameter at least four, we use Theorem~\ref{teo1}. The following result, which is based on the results of~\cite{TCDV10}, summarizes what is known for trees with diameter three.
\begin{Lema}\label{teodiam3}\cite{TCDV10}
Consider an $n$-vertex tree $T(a,b)$, where $a \geq b \geq 1$ and $n=a+b+2$.
\begin{itemize}
\item[(a)] The characteristic polynomial of the Laplacian matrix of $T(a,b)$ is given by
$$P_{T(a,b)}(x)=p_{a,b}(x)\cdot(x-1)^{n-4}\cdot x$$ where $p_{a,b}(x)=x^3-(b+a+4)x^2+(ab+2b+2a+5)x-b-a-2.$

\vspace{4pt}

\item[(b)] The tree $T(a,b)$ has $\sigma=2$ Laplacian eigenvalues that are larger than the average degree $2-2/n$, and $S_2(T(a,b))=y_1+y_2=n+2-y_3$, where $y_1 \geq y_2 \geq y_3$ are the roots of $p_{a,b}(x)$.

\vspace{4pt}

\item[(c)] The algebraic connectivity of $T(a,b)$ is equal to $y_3>2/n$ and, for $k \geq 2$,
$$S_k(T(a,b)) < n+k-\frac{2}{n}.$$

\end{itemize}
\end{Lema}

\begin{proof}[Proof of Theorem~\ref{order}] For $n\geq 6$, consider an $n$-vertex tree $T(a,b)$ with diameter three, where $a \geq b \geq 1$, and an $n$-vertex tree $T$ with diameter at least four. Let $\sigma$ denote the number of eigenvalues of $T$ that are larger than the average $\overline{d}$. By Proposition~\ref{prop_1} and Lemma~\ref{teodiam3}, the difference between the Laplacian energies of $T(a,b)$ and $T$ is given by
    \begin{eqnarray}
    LE(T(a,b))-LE(T)&=&2(n+2-y_3)-4\overline{d}-2S_\sigma(T)+2\sigma\overline{d} \nonumber\\
    &\geq&2n+4-2y_3-8+\frac{8}{n}-2\left(n-2+2\sigma-\frac{2\sigma}{n}\right)+4\sigma-\frac{4\sigma}{n} \label{mona1}\\
    &=&-2y_3+\frac{8}{n} \nonumber.
    \end{eqnarray}
Observe the use of Theorem~\ref{teo1} to obtain~\eqref{mona1}. In particular, in order to show that $LE(T(a,b))-LE(T)\geq0$, it suffices to prove that $y_3\leq\frac{4}{n}$. Note that $p_{a,b}(x)<0$ for all $x<y_3$, so that it suffices to show that $p_{a,b}(4/n) \geq 0$. Let $a=\left\lceil\frac{n-2}{2}\right\rceil+k$ and $b=\left\lfloor\frac{n-2}{2}\right\rfloor-k$, where $k \in \{0,\ldots,\left\lfloor\frac{n-2}{2}\right\rfloor-1\}$.

Observe that, for $n$ even, $$p_{a,b}\left(\frac{4}{n}\right)=\frac{4(n^3-k^2n^2-2n^2-8n+16)}{n^3},$$ which is nonnegative, for nonnegative integers $k$, if and only if
$$k\leq \left\lfloor \sqrt{n-2-\frac{8}{n}+\frac{16}{n^2}} \right\rfloor.$$
To prove the theorem in this case, we need to show that $\left\lfloor \sqrt{n-2-\frac{8}{n}+\frac{16}{n^2}} \right\rfloor =  \left\lfloor \sqrt{n-3} \right\rfloor$, which follows from a few algebraic manipulations that are included in the appendix. This concludes the proof of Theorem~\ref{order} for $n$ even, since Theorem~\ref{orderTrevisan} tells us that the Laplacian energy of $T\left(\left\lceil\frac{n-2}{2}\right\rceil+k,\left\lfloor\frac{n-2}{2}\right\rfloor-k\right)$ decreases as $k \in \{0,\ldots,\left\lfloor \sqrt{n-3} \right\rfloor-1\}$ increases.

If $n$ is odd, we have $a=(n-1)/2+k$ and $b=(n-3)/2-k$, so that
$$p_{a,b}\left(\frac{4}{n}\right)=\frac{4(n^3-k^2n^2-n^2k-8n+16)-9n^2}{n^3},$$
which is nonnegative if and only if
$$k\leq\left\lfloor\sqrt{n-2-\frac{8}{n}+\frac{16}{n^2}}-\frac{1}{2} \right\rfloor.$$
The result follows from the fact that
$\left\lfloor\sqrt{n-2-\frac{8}{n}+\frac{16}{n^2}}-\frac{1}{2} \right\rfloor=\left\lfloor\sqrt{n-3}+\frac{1}{2}\right\rfloor$, if $n=p^2-p+3$ for some integer $p\geq3$, and that $\left\lfloor\sqrt{n-2-\frac{8}{n}+\frac{16}{n^2}}-\frac{1}{2} \right\rfloor=\left\lfloor\sqrt{n-3}-\frac{1}{2}\right\rfloor$ for all other odd values of $n \geq 7$.
 \end{proof}

\begin{Example}
For $n=42$, the star $\mathcal{S}_{42}$ has Laplacian energy
(truncated to four decimal places) equal to 80.0952. The other
$f(42)=1+\lfloor\sqrt{42-3}\rfloor=7$ trees with largest Laplacian
energy are given by Table \ref{tabelaEx}. Theorem~\ref{teo1}
implies that the Laplacian energy of any tree with diameter larger
than or equal 4 with 42 vertices is less than 80. The tree
$T(27,13)$, which has the largest Laplacian energy among all
remaining trees with diameter three, has energy 79.9959. This is
less than the given upper bound.
\begin{table}[h]\label{tabelaEx}    \begin{center} \begin{tabular}{|c|c|c|c|}  \hline
    ranking & tree & $a-b$ & $LE$\\ \hline
    2 & $T(20,20)$ & 0 & 80.0159\\ \hline
    3 & $T(21,19)$ & 2 & 80.0155\\ \hline
    4 & $T(22,18)$ & 4 & 80.0144\\ \hline
    5 & $T(23,17)$ & 6 & 80.0125\\ \hline
    6 & $T(24,16)$ & 8 & 80.0098\\ \hline
    7 & $T(25,15)$ & 10 & 80.0062\\ \hline
    8 & $T(26,14)$ & 12 & 80.0016\\ \hline
    \end{tabular} \caption{Example for $n=42$.} \end{center} \end{table}

For $n=43$, in addition to the star $\mathcal{S}_{43}$, we can identify the next $f(43)=\lfloor\sqrt{43-3}+\frac{1}{2}\rfloor=6$ trees with highest Laplacian energy. For $n=45$, we have $f(45)=\lfloor\sqrt{45-3}+\frac{3}{2}\rfloor=7$, because 45 is equal to $p^2-p+3$ for $p=7$.

\end{Example}

In light of Theorem~\ref{order}, it is natural to ask whether the Laplacian energy of any tree with diameter three is larger than the Laplacian energy of all trees with the same number of vertices and larger diameter. As a matter of fact, this holds for all trees with up to 15 vertices (see the database on the website~\cite{siteC}). However, this is not true in general: for any fixed $n \geq 16$, there is a tree with $n$ vertices and diameter 4 whose Laplacian energy is larger than the Laplacian energy of $T(n-3, 1)$, the tree with the smallest Laplacian energy among all $n$-vertex trees with diameter three. Given $n\geq 16$ and $k = \lfloor \frac{n}{3} \rfloor$, we consider the tree $F(n,k)$ (see Figure~\ref{badtree}), where a central vertex is incident with $n-2k-1$ leaves (which we call pendants), and with two vertices of degree $k$, each of which is incident with $k-1$ pendants.
\begin{figure}[h!]\label{badtree}
    \begin{center}
        \includegraphics[width=0.25\linewidth]{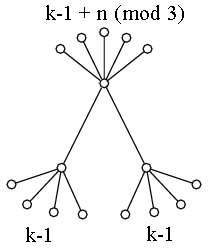}
        \caption{Tree $F(n,k), k = \lfloor \frac{n}{3} \rfloor$}
        \label{badtree}
    \end{center}
\end{figure}\\

\begin{table}[h!]\label{tableasmall}
\begin{center} \begin{tabular}{|c|c|c|}  \hline
    $(n,k)$ & $LE(F(n,k))$ & $LE(T(n-3,1))$ \\ \hline
    (16,5) & 27.6803 & 27.6739 \\ \hline
    (17,5) & 29.6830 & 29.6485 \\ \hline
    (18,6) & 31.6815 & 31.6259 \\ \hline
    (21,7) & 37.6982 & 37.5709\\ \hline
    \end{tabular} \caption{Small values of ($n,k$).} \end{center}
\end{table}

For ($n, k$) = (16, 5), (17, 5), (18, 6) or (21, 7), we show in
Table \ref{tableasmall} (whose computation may be done by hand)
that $LE(F(n,k)) > LE(T(n-3,1))$.

\begin{Teo}\label{thm_counter}
 Let $n\geq 16$ and $k =\lfloor \frac{n}{3} \rfloor$. Then $LE(F(n,k)) > LE(T(n-3,1))$.
\end{Teo}
We may prove that the same relationship holds in general. The proof of this result relies on careful calculations to find appropriate upper and lower bounds on the roots of the characteristic polynomials of the Laplacian matrices of $F(n,k)$ and $T(n-3,1)$, and we include it as an appendix.

To conclude this section, we observe that Theorem~\ref{teo1} leads to a new upper bound on $S_k(G)$ for a general graph $G$. Wang~\etal~\cite[Theorem~3.1]{Wang} have shown that, if $G$ is a connected graph with $n$ vertices and $1 \leq k \leq n$, we have
$$S_k(G) \leq 2e(G)+2k-n-\frac{2k-2}{n}.$$
With essentially the same arguments, but making use of our improved bound on $S_k(T)$, we may obtain the following result. For completeness, we include a proof of this result.
\begin{Teo}\label{teo_cyclic}
For a graph $G$ with $n\geq 6$ vertices and diameter $d \geq 4$, we have
$$S_k(G) < 2e(G)-n+2k-\frac{2k}{n}.$$
\end{Teo}

 The main technical tool is the following classical result by Wielandt~\cite{W55}.
\begin{Teo}{\label{teoautmat}}
Let $A, B, C$ be hermitian matrices of order $n$ such that $A=B+C$. Then
$$\sum_{i \in I}\lambda_i(A) \leq \sum_{i=1}^{|I|}\lambda_i(B) + \sum_{i \in I}\lambda_i(C)$$
for any subset $I\subset\{1,2,...,n\}$.
\end{Teo}

\begin{proof}[Proof of Theorem~\ref{teo_cyclic}]
Clearly, it suffices to show this result for connected graphs and, given an $n$-vertex graph $G$, we let $c=c(G)$ be the nonegative integer such that $|E(G)|=n-1+c$. With this, the upper bound may be rewritten as
$$S_k(G) < n-2+2(k+c)-\frac{2k}{n}.$$
Our proof is by induction on $c$. If $c=0$, then $G$ is a tree and the result is simply Theorem~\ref{teo1}.

Assume that $G$ is an $n$-vertex graph with $c \geq 1$, so that $G$ contains a cycle. Let $e$ be an edge in a cycle of $G$ and consider the connected graph $G'=G-e$. Clearly, the diameter of $G'$ is at least the diameter of $G$. With an appropriate ordering of the vertices of $G'$, we may decompose the Laplacian matrix of $G$ as
\begin{equation*}
L(G)=L(G')+M=L(G')+\begin{pmatrix}
    M_{*} & 0 \\ 0 & 0
    \end{pmatrix},
\end{equation*}
where $M_{*}=\begin{pmatrix} 1 & -1 \\ -1 & 1 \end{pmatrix}$ has eigenvalues $0$ and $2$. Theorem~\ref{teoautmat} implies that
\begin{eqnarray*}
S_k(G) \leq S_k(G') + S_k(M) \leq S_k(G') + 2,
\end{eqnarray*}
and our result follows by induction, since
\begin{eqnarray*}
S_k(G') + 2 <  n-2+2(k+c-1)-\frac{2k}{n} + 2 = n-2+2(k+c)-\frac{2k}{n},
\end{eqnarray*}
which is precisely the required upper bound.
\end{proof}

Using the relation $e(G)=n-1+c$, Theorem~\ref{teo_brouwer} follows directly from the above result if we observe that $$2k-\frac{2k}{n} + c-1 <  {k+1 \choose 2}$$
for every $k$ satisfying
$$\frac{3-\frac{4}{n}+\sqrt{(3-\frac{4}{n})^2+8c-8}}{2} \leq k \leq n.$$

To conclude this section, we justify our claim that Theorem~\ref{teo_brouwer} is an improvement on Corollary~3.2 in~\cite{Wang} by showing that
$$\left\lfloor\frac{3n-4+\sqrt{9n^2-24n+16+8e(G)n^2-8n^3}}{2n} \right\rfloor<\left\lfloor\frac{3n-4+\sqrt{9n^2-8n+16+8e(G)n^2-8n^3}}{2n} \right\rfloor$$
whenever $n=11+j$ and $e(G)=(6+j)(7+j)/2$. To this end, observe that the expressions inside the square root on the left-hand side and on the right-hand side of this inequality are equal to
$$f(j)=10521+8224 j+2424 {j}^{2}+320 {j}^{3}+16\,{j}^{4} \textrm{ and }g(j)=10697+8256 j+2424{j}^{2}+320{j}^{3}+16{j}^{4},$$
respectively. Note that the coefficients of $f(j)$ are smaller than or equal to the corresponding coefficients of $g(j)$. Also comparing coefficients, we see that
$$h(j)=(4j^{2}+40j+103)^{2}=16 j^{4}+320 j^{3}+2424 j^{2}+8240 j+10609$$
satisfies $f(j)<h(j)<g(j)$ for every $j \geq 0$. Our claim now follows from the fact that
$\frac{3n-4+\sqrt{h(j)}}{2n}=j+6$, and hence  $\left\lfloor \frac{3n-4+\sqrt{f(j)}}{2n} \right\rfloor \leq j+5$ and  $\left\lfloor \frac{3n-4+\sqrt{g(j)}}{2n} \right\rfloor \geq j+6$. In fact, these two inequalities hold with equality for every nonnegative integer $j$.

\section{Proof of Theorem~\ref{teo1} - Overview of the proof}\label{sec_overview}

In order to describe the structure of the proof of Theorem~\ref{teo1}, we compare it with the proof of Theorem~1.1 in~\cite{FHRT11}, where we showed that every $n$-vertex tree satisfies the inequality~\eqref{eq_bound}. This allows us to point out similarities and differences in the proofs of these two results.

The first observation is that, to obtain a general bound on the sum of the largest Laplacian eigenvalues of a tree $T$, it suffices to verify that it holds in the case $k=\sigma$, where $\sigma$ is the number of Laplacian eigenvalues of $T$ that are larger than the average degree $\overline{d}$. For future reference, we state this result in terms of the bound in the current paper. We omit the proof, which is based on straightforward calculations and follows the same steps as~\cite[Lemma~3.1]{FHRT11}.
\begin{Lema}{\label{lemassigma}} Let $T$ be an $n$-vertex tree with $\sigma$ eigenvalues larger than the average degree $\overline{d}=2-2/n$. If the inequality $$S_k(T) < n-2+2k-\frac{2k}{n}$$ holds for $k=\sigma$, then it holds for every $k \in \{1,\ldots,n\}$.
\end{Lema}

In~\cite{FHRT11} the bound~\eqref{eq_bound} was obtained in two steps, which are inspired by the work of Haemers~\etal~\cite{HMT10}. Given a tree $T$, the basic idea is to apply induction by removing an edge of $T$, so as to obtain the bound for its Laplacian energy in terms of the bound for two of its subtrees. As it turns out, this general argument suffices for the results in~\cite{HMT10}, and, for the upper bound in~\eqref{eq_bound}, it may be carried out directly whenever we can choose an edge $e$ for which the two components of $T-e$ have diameter at least three. The remaining trees, namely all trees such that the removal of any edge produces at least one star, were considered separately. In some sense, these are the trees for which the bound is tightest, and the basic tool to bound their Laplacian energy was to find information about $\sigma$ in order to estimate $S_\sigma(T)$.

In this paper, the structure of the proof is along the same line. Instead of using induction, we deal with the general case by contradiction, considering a minimum counterexample. Once again, the trees that have been considered separately in~\cite{FHRT11} do not fit the general argument and are analyzed with different arguments, which either rely on computing characteristic polynomials directly, and on estimating their roots, or depend on the decomposition of the Laplacian matrix into simpler matrices.

To compute the characteristic polynomials of the Laplacian matrix, we use a straightforward translation of an algorithm due to Jacobs, Machado and Trevisan~\cite{JMT05} to the context of Laplacian matrices (see~\cite{FHRT11} for details about this adaptation; the algorithm in~\cite{JMT05} calculates the characteristic polynomial of the adjacency matrix). One of the advantages of this method is that it can be slightly modified in such a way that we can determine, for any $\alpha \in \mathbb{R}$, the number of Laplacian eigenvalues of a tree $T$ that are larger than $\alpha$, equal to $\alpha$ and smaller than $\alpha$, respectively (see~\cite{JT11} for the original work, which deals with eigenvalues of the adjacency matrix, and~\cite{FHRT11} for the adaptation to the Laplacian scenario). The algorithm is initiated as follows: the vertices of the input tree $T$ are labelled $1,\ldots,|V(T)|$ (the vertex with largest label is called the root) so that the sequence of labels of the vertices on any path between a leaf and the root is increasing (in other words, the vertices are ordered bottom-up with respect to the root). The quantity $a(v)=d(v)-\alpha$ is associated with each vertex $v \in V(T)$, where $d(v)$ is the degree of $v$ in $T$. The vertices are then processed one by one according to this ordering, starting with the vertex with label $1$. Nothing is done if $v$ is has no children, otherwise $a(v)$ is updated to $a(v)-\sum \frac{1}{a(c)}$, where the sum runs over all children of $v$. In case $a(c)=0$ for some child $c$ of $v$, the algorithm sets $a(v)=-1/2$ and $a(c)=2$, and it removes the edge between $v$ and its parent if $v$ is not the root.  After all vertices are processed, the number of vertices for which $a(v)>0$ gives the number of Laplacian eigenvalues of $T$ that are larger than $\alpha$, while the number of vertices for which $a(v)<0$ gives the number of Laplacian eigenvalues of $T$ that are smaller than $\alpha$. Hence the number of eigenvalues that are equal to $\alpha$ is the number of vertices such that $a(v)=0$.

To estimate the sum of the eigenvalues of the Laplacian matrix of a tree in terms of simpler matrices, we need to introduce the special family of trees under consideration.  In~\cite{FHRT11} they were called $S \& S$-trees (Star-and-Star trees) or $SNS$-trees (Star-NonStar trees) according to whether they could be split into two stars by the removal of an edge that is not incident with a leaf, or whether the removal of any such edge would produce two components, one of which is a star and one of which has diameter at least three (i.e., a `non-star'). In the current paper, however, we shall consider these two families of trees together, as the set of $S\& S$-trees with diameter at least four has a very particular structure: if the tree is obtained from two stars by adding an edge between the center of one of the stars and a leaf of the other, it is a \emph{double broom with diameter four}, that is, a tree with diameter four consisting of a path on five vertices whose central vertex has degree two, while the other two nonleaf vertices may have arbitrary degree. On the other hand, if a tree is obtained from two stars by adding an edge joining leaves, it is a \emph{double broom with diameter five}, which is a tree with diameter five consisting of a path on six vertices such that an arbitrary number of leaves may be appended to the second and the fifth vertices.

Following the notation and terminology of~\cite{FHRT11}, it is easy to see that any tree in the special family $\mathfrak{F}$, which comprises $SNS$-trees and $S\&S$-trees with diameter at least four, may be viewed as a tree with a \emph{root vertex} $v_0$ with which three different types of branches may be incident:
\begin{itemize}

\item[(i)] a single vertex, also called a \emph{pendant}, which we call a branch of \emph{type 0};

\item[(ii)] a tree of height one, called a branch of \emph{type 1};

\item[(iii)] a tree of height two whose root has degree one, called a branch of \emph{type 2}.

\end{itemize}
Moreover, the combined number of branches of type 1 and 2 in a tree in $\mathfrak{F}$ is at least two.

It follows immediately from the definition that every tree in $\mathfrak{F}$ has one of three possible diameters: the absence of branches of type 2 implies diameter four; exactly one branch of type 2 leads to diameter five; two or more branches of type 2 give diameter six. Based on this, we say that a tree in $\mathfrak{F}$ lies in $\mathfrak{F}_4$, $\mathfrak{F}_5$ or $\mathfrak{F}_6$ according to its diameter.

\begin{figure}[h,t]
  \begin{center}
  \includegraphics[width=2.5 in, height=1.5 in]{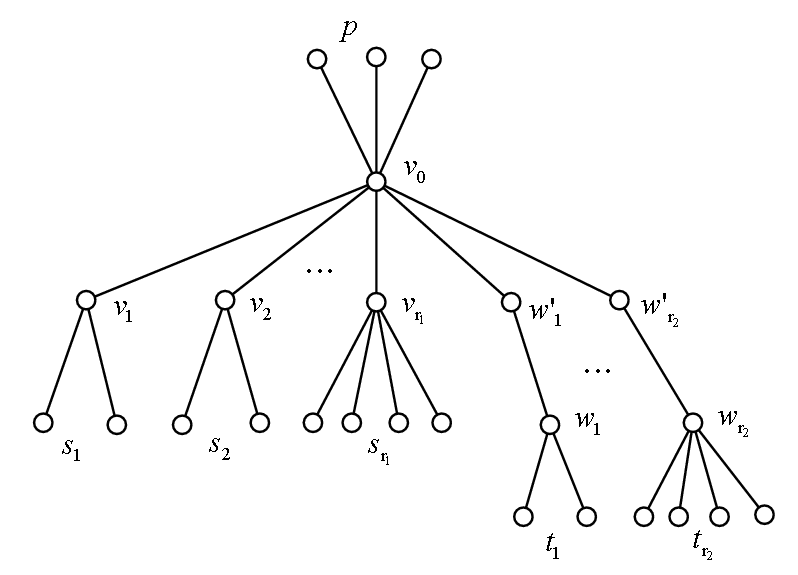}
  \caption{Trees in $\mathfrak{F}$}\label{prototype}
\end{center}\end{figure}

Henceforth, when referring to a tree $T \in \mathfrak{F}$, we use the following notation,  which is depicted in Figure~\ref{prototype}. The tree $T$ has a central node $v_0$ that is adjacent to three types of branches: $p\geq 0$ pendants; $r_1\geq 0$ branches of type 1, which are rooted at vertices $v_1, \ldots, v_{r_1}$, where the branch rooted at $v_i$ has $s_i \geq 1$ leaves; and
$r_2 \geq 0$ branches of type 2, such that the $j$-th branch is rooted at a vertex $\ww_j$, whose single neighbor $\w_j$ in the branch is adjacent to
$t_j \geq 1$ leaves. We shall suppose that $s_1\leq \cdots \leq s_{r_1}$ and $t_1\leq \cdots \leq t_{r_2}$. Clearly, we have $r_1+r_2 \geq 2$, and the total number vertices of a particular tree is $n=p+r_1+2r_2+1+\sum_{i=1}^{r_1} s_i + \sum_{j=1}^{r_2} t_j$.

The three types of branches under consideration are special instances of \emph{generalized Bethe trees}, which are rooted trees such that vertices at the same distance from the root have equal degrees. Rojo~\cite{Roj08} has characterized the Laplacian spectrum of generalized Bethe trees whose roots are adjacent to a common root. Since the branches of a tree in $\mathfrak{F}$ are all connected to a common root $v_0$, his result can be applied to our framework .

Let $T$ be a tree in $\mathfrak{F}$. Adapting the notation of~\cite{Roj08}, we define, for $i=1,2,...,r_1$, the matrix $T_i$ of order 2 associated with the $i$-th branch of type 1:
\begin{equation}{\label{eqndiam4ti}}
    T_i = \begin{pmatrix}
    1 & \sqrt{s_i} \\ \sqrt{s_i} & s_i + 1
    \end{pmatrix}.
    \end{equation}
Moreover, for $j=1,2,...,r_2$, we define the matrix $Q_i$ of order 3 associated with the $j$-th branch of type 2:
\begin{equation} \label{eqndiam4qj}
Q_j = \begin{pmatrix}
    1 & \sqrt{t_j} & 0 \\ \sqrt{t_j} & t_j+1 & 1 \\ 0 & 1 & 2
    \end{pmatrix}.
\end{equation}

Let $M$ be the matrix defined as one of the matrices $M_p$ or $M_{\bar{p}}$ below, according to whether $T$ has pendant vertices (matrix $M_p$) or not (matrix $M_{\bar{p}}$):
\begin{eqnarray}
    M_p = \begin{pmatrix}
    1 &&&&&&& \sqrt{p} \\
    &T_1 &&&&&& u_1 \\
    && \ddots &&&&& \vdots \\
    &&& T_{r_1} & &&&u_{r_1} \\
    &&&& Q_1 &&& \overline{u}_{1}\\
    &&&&& \ddots && \vdots \\
    &&&&&& Q_{r_2} &\overline{u}_{r_2}\\
    \sqrt{p} & u_1^T & \dots & u_{r_1}^T & \overline{u}_{1}^T & \ldots & \overline{u}_{r_2}^T & \delta,
    \end{pmatrix} \label{Tezona0}
\end{eqnarray}
\begin{eqnarray}
    M_{\bar{p}} = \begin{pmatrix}
    T_1 &&&&&& u_1 \\
    & \ddots &&&&& \vdots \\
    && T_{r_1} & &&&u_{r_1} \\
    &&& Q_1 &&& \overline{u}_{1}\\
    &&&& \ddots && \vdots \\
    &&&&& Q_{r_2} &\overline{u}_{r_2}\\
    u_1^T & \dots & u_{r_1}^T & \overline{u}_{1}^T & \ldots & \overline{u}_{r_2}^T & \delta,
    \end{pmatrix} \label{Tezona1}
    \end{eqnarray}
where $T_i$, $i=1,\ldots, r_1$, and $Q_j$, $j=1,\ldots, r_2$, are defined by equations (\ref{eqndiam4ti})
and (\ref{eqndiam4qj}), respectively, $u_i=[0, 1]^T$, $\overline{u}_j=[0, 0, 1]^T$ and $\delta=r_1+r_2+p$ is the degree of $v_0$ in $T$.

The following result, which explicits the connection between the Laplacian spectrum of $T$ and the spectrum of the matrix $M$, can be read from Theorem 2 in~\cite{Roj08}.
\begin{Lema}[\cite{Roj08}]\label{decomp_Rojo} The Laplacian spectrum of $T \in \mathfrak{F}$ is the multiset given by the union of the spectrum of the matrix $M$ defined in equation (\ref{Tezona0}) or (\ref{Tezona1}) and of a multiset where all the elements are equal to 1.
\end{Lema}

One of the implications of this result is that all the Laplacian eigenvalues that are larger than their average $\overline{d}$ are also eigenvalues of the matrix $M$. Hence, in order to estimate the sum of the largest eigenvalues of $T$, it suffices to estimate the eigenvalues of $M$. To this end, we decompose $M$ as a sum of two well-understood matrices, to which we may apply Theorem~\ref{teoautmat}.

If $T$ contains $p \geq 1$ pendant vertices (hence $M=M_p$), we let

\begin{eqnarray}\label{eqndiamgeralcp}
M_p &=&  A_p + B_p \\
&=& \begin{pmatrix}
    1 &&&&&&&  \\
    &T_1 &&&&&&  \\
    && \ddots &&&&& \\
    &&& T_{r_1} & &&& \\
    &&&& Q_1 &&& \\
    &&&&& \ddots &&  \\
    &&&&&& Q_{r_2} &\\
    &&&&&&& \delta
    \end{pmatrix}  +
\begin{pmatrix}
    &&&&&&& \sqrt{p} \\
    &&&&&&& u_1 \\
    &&&&&&& \vdots \\
    &&& & &&&u_{r_1} \\
    &&&&&&& \overline{u}_{1}\\
    &&&&&&& \vdots \\
    &&&&&&&\overline{u}_{r_2}\\
    \sqrt{p} & u_1^T & \dots & u_{r_1}^T & \overline{u}_{1}^T & \ldots & \overline{u}_{r_2}^T & 0
    \end{pmatrix}. \nonumber
    \end{eqnarray}
Analogously, if there are no pendant vertices in $T$ (hence $M=M_{\bar{p}}$), we let
  \begin{eqnarray}{\label{eqndiamgeralsp}}
    M_{\bar{p}} &=&  A_{\bar{p}}+B_{\bar{p}}\\
    &=&
\begin{pmatrix}
    T_1 &&&&&& \\
    & \ddots &&&&&\\
    && T_{r_1} & &&&\\
    &&& Q_1 &&&\\
    &&&& \ddots && \\
    &&&&& Q_{r_2}& \\
    &&&&&&\delta
    \end{pmatrix}
+  \begin{pmatrix}
    &&&&&& u_1 \\
    &&&&&& \vdots \\
    && & &&&u_{r_1} \\
    &&& &&& \overline{u}_{1}\\
    &&&& && \vdots \\
    &&&&& &\overline{u}_{r_2}\\
    u_1^T & \dots & u_{r_1}^T & \overline{u}_{1}^T & \ldots & \overline{u}_{r_2}^T & 0
    \end{pmatrix}.
\nonumber
    \end{eqnarray}

In order to apply Theorem~\ref{teoautmat} to the above decompositions of $M$, we need the spectra of the matrices $A$ and $B$ (which have subindices ``$p$'' or ``$\bar{p}$'' according to whether $M$ is equal to $M_p$ or $M_{\bar{p}}$). Let  $(\sqrt{b_1},\sqrt{b_2},\dots,\sqrt{b_{n-1}},0)$ denote the last row (and last column) of $B$. It is easy to prove (see Lemma~4.2 in~\cite{FHRT11}) that the spectrum of $B$ is given by $\{\sqrt{\delta},0,0,\ldots,0,-\sqrt{\delta}\}$, where $\delta=\sum_{i=1}^{n-1}b_i$. This allows us to relate the largest Laplacian eigenvalues of $M$ with the largest eigenvalues of $A$ in the decompositions (\ref{eqndiamgeralcp}) and (\ref{eqndiamgeralsp}). For future reference, this is stated as the following result, which is obtained by applying Theorem~\ref{teoautmat} to $A+B$ with index set $I=\{1,\ldots,k,|M|\}$, where $|M|$ denotes the order of $M$.
\begin{Lema}\label{aplic_decomp_legal}
Let $T \in \mathfrak{F}$ be an $n$-vertex tree with Laplacian eigenvalues $\mu_1 \geq \cdots \geq  \mu_n$, where exactly $\sigma$ of them are larger than their average $\overline{d}$. Let $M=A+B$ be the matrix decomposition in (\ref{eqndiamgeralcp}) or (\ref{eqndiamgeralsp}) associated with it. Then, for $1 \leq k \leq \sigma$, we have
\begin{eqnarray}{\label{eqndiamgeral}}
 S_{k}(T)=\sum_{i=1}^{k} \mu_k \leq S_{k+1}(A)=\sum_{i=1}^{k+1} \lambda_i(A).
\end{eqnarray}
\end{Lema}

The spectrum of $A$ may be read directly from the spectra of its submatrices $T_i$ and $Q_j$, which we now describe.
\begin{Lema}{\label{lemamatriza2}}~\cite[Lemma 4.4]{FHRT11} For every integer $ s\geq 1$, the matrix
\begin{equation*}
T = \begin{pmatrix}
1 & \sqrt{s} \\ \sqrt{s} & s+1
\end{pmatrix}
\end{equation*}
has eigenvalues $x_1 > x_2$ satisfying
\begin{eqnarray*}
2 < &x_1& < 2+s-\frac{1}{2+s},\\
0 < &x_2& < \frac{1}{2}.
\end{eqnarray*}
\end{Lema}

\begin{Lema}{\label{lemamatriza3}}~\cite[Refinement of Lemma~4.5]{FHRT11} For every integer $t \geq 1$, the eigenvalues $y_1 > y_2 > y_3$ of the matrix
 \begin{equation*}
 Q = \begin{pmatrix}
 1 & \sqrt{t} & 0 \\ \sqrt{t} & t+1 & 1 \\ 0 & 1 & 2
 \end{pmatrix}
 \end{equation*}
 satisfy the following:
 \begin{eqnarray*}
 2 < &y_1& < t+2+\frac{1}{4t} \\
 \frac32< &y_2& < 2 \\
 1> &y_3& > \begin{cases}
    0.19 &  \textrm{ for } t=1 \\
    \frac{1}{4t} & \textrm{ for } t\geq 2 \end{cases}  \\
 y_1+&y_2&=t+4-y_3
 \end{eqnarray*}
Furthermore, the second eigenvalue increases as $t$ increases.
 \end{Lema}

\section{Proof of Theorem~\ref{teo1} - Particular cases}\label{sec_particular}

In this section we prove Theorem~\ref{teo1} for the particular family of trees $\mathfrak{F}$ defined in the previous section.

\subsection{Trees in $\mathfrak{F}_4$}

Following the notation used in Figure~\ref{prototype}, a tree $T$ in $\mathfrak{F}_4$ has a central node $v_0$, which is adjacent to $p\geq 0$ pendants and $r_1\geq 2$ branches of type 1, which are rooted at vertices $v_1, \ldots, v_{r_1}$, where the branch rooted at $v_i$ has $s_i \geq 1$ leaves. There are no branches of type 2.

The following result summarizes some spectral properties of the trees in $\mathfrak{F}_4$ obtained in~\cite{FHRT11}.
\begin{Lema}{\label{lemadiam4}}\cite{FHRT11}
A tree in $\mathfrak{F}_4$ with $r\geq2$ branches of type 1 has $r$ or $r+1$ eigenvalues larger than $\overline{d}$. Moreover, if the root $v$ of such a tree $T$ is incident with a pendant, or if $v$ is incident with at least three branches of type 1 with two or more leaves, then $\sigma=r+1$.
\end{Lema}

\begin{proof}[Proof of Theorem~\ref{teo1} for $\mathfrak{F}_4$] To prove Theorem~\ref{teo1} for $\mathfrak{F}_4$, we shall consider two basic cases, according to the value of $\sigma(T)$, which may be equal to $r$ or $r+1$ because of Lemma~\ref{lemadiam4}. By Lemma \ref{lemassigma}, it suffices to check the validity of~\eqref{new_bound} for $k=\sigma$.

\vspace{5pt}

\noindent \emph{Case 1 $(\sigma=r+1)$:} Let $T$ be a tree in $\mathfrak{F}_4$ such that $\sigma(T)=r+1$. The argument used to show that~\eqref{new_bound}  holds depends on whether $r \geq 3$ or $r=2$.

\vspace{5pt}

\noindent \emph{Case 1.1 $(r \geq 3)$:} We shall apply Lemma~\ref{aplic_decomp_legal} to $T$ with $k=r+1$. It implies that
$$S_{r+1}(T) \leq S_{r+2}(A),$$
where $A$ is the matrix defined in (\ref{eqndiamgeralcp}) or (\ref{eqndiamgeralsp}) according to whether $v_0$ is adjacent to pendantss (i.e. $A=A_p$) or not (i.e. $A=A_{\bar{p}}$). We shall use the following inequality, which may be proved easily (for instance with elementary calculus).
\begin{Lema}{\label{lemacotasomainvsi2}} For $a_i\geq1$, $1\leq i\leq r$, such that $\sum_{i=1}^r a_i=c$, we have $$\sum_{i=1}^{r}\frac{1}{a_i+2}\geq\frac{r^2}{c+2r}.$$
\end{Lema}

By Lemma~\ref{lemamatriza2} we know that the $r+1$ largest eigenvalues of $A_{\bar{p}}$ are the degree $\delta$ of $v_0$ and the largest eigenvalues of $T_1, \ldots,T_r$. Finally, $\lambda_{r+2}(A_{\bar{p}})$ is the largest among all smallest eigenvalues of $T_1, \ldots,T_r$, which is equal to $\frac{1}{2}\left(s+2-\sqrt{s^2+4s} \right)$ for some $s \geq 1$ and is therefore smaller than $1/2$. On the other hand, the $r+2$ largest eigenvalues of $A_p$ are, in nonincreasing order, $\delta$, the largest eigenvalues of $T_1,\ldots,T_r$, and 1. As a consequence, we know that, in both cases, $S_{r+1}(T)$ satisfies
\begin{eqnarray}
    S_{r+1} (T) &\leq& \delta + \sum_{i=1}^{r}\left(\frac{s_i+2+\sqrt{s_i^2+4s_i}}{2}\right) + 1 \label{mona7} \\
    &<& \delta + 1 + \sum_{i=1}^{r}\left( s_i+2-\frac{1}{s_i+2}\right)  \label{eqaux1}\\
    &=& \delta + 1 + 2r + \sum_{i=1}^{r}s_i -\sum_{i=1}^{r}\frac{1}{s_i+2} \nonumber \\
    &=& n + 2r -\sum_{i=1}^{r}\frac{1}{s_i+2} \label{eqaux2}\\
    &\leq& n-1 + 2(r+1)-1 -\frac{r^2}{\sum_{i=1}^{r} s_i + 2r} \label{eqaux3}\\
    &=& n-1 +2(r+1)-1 -\frac{r^2}{n+r-p-1}. \label{demdiam4compfinal}
    \end{eqnarray}
For~(\ref{eqaux1}), we used Lemma~\ref{lemamatriza2}, for~(\ref{eqaux2}) and~(\ref{demdiam4compfinal}), we used the relation $n=1+p+r+\sum_{i=1}^r s_i$, while~(\ref{eqaux3}) comes from Lemma~\ref{lemacotasomainvsi2}. We shall obtain our result if we prove that
\begin{equation}\label{cooper1}\frac{r^2}{n+r-p-1} \geq \frac{2(r+1)}{n}.\end{equation}
Assuming that $p \geq 1$, since $r \geq 3$, we have
\begin{eqnarray}
&\frac{r^2}{n+r-p-1} \geq \frac{2r+2}{n} \label{mona8}\\
\Longleftrightarrow &r^2n\geq2rn+2r^2-2rp+2n-2p-2 \nonumber\\
\Longleftrightarrow &n(r^2-2r-2)\geq2(r^2-rp-p-1) \nonumber\\
\Longleftarrow &n(r^2-2r-2)\geq2(r^2-r-2) \nonumber\\
\Longleftrightarrow &n\geq2\frac{r^2-r-2}{r^2-2r-2}=2+\frac{2r}{r^2-2r-2} \nonumber\\
\Longleftarrow &n\geq2+2\cdot3=8\nonumber,
\end{eqnarray}
which holds as $n=1+p+r+\sum s_i\geq1+p+2r \geq8$.

In the case $p=0$, we may replace the additive term 1 by 1/2 in the upper bound given in~\eqref{mona7}. In particular, instead of \eqref{mona8}, it suffices to show that
\begin{eqnarray*}
&\frac{r^2}{n+r-p-1} +\frac{1}{2} \geq \frac{2r+2}{n}\\
 \Longleftrightarrow &\frac{2r+2}{n}-\frac{r^2}{n+r-p-1} \leq \frac{1}{2} \\
 \Longleftrightarrow &\frac{n(2+2r-r^2)+2r^2}{n(n+r-1)} \leq \frac{1}{2}.
\end{eqnarray*}
We know that $2+2r-r^2<0$, as $r \geq 3$ and that $r<n/2$. It follows that
$$\frac{n(2+2r-r^2)+2r^2}{n(n+r-1)}\leq \frac{2r^2}{n(n+2)} \leq \frac{1}{2},$$
as required.

\vspace{5pt}

\noindent \emph{Case 1.2 $(r = 2)$:} We assume that $T$ is an $n$-vertex tree such that $\sigma=r+1$ and $r=2$. Here, the result does not follow from~\eqref{cooper1}. We shall address this case directly through the characteristic polynomial.

Denote $s_1=a\geq s_2=b\geq 1$ the number of leaves in each of the two branches of type 1 and $p \geq 0$ the number of pendants. We have $n=p+a+b+3$ vertices. The characteristic polynomial of the Laplacian matrix of this tree is given by
\begin{equation*}
P_T(x)=q_{p,a,b}(x)\cdot(x-1)^{n-6}\cdot x,
\end{equation*}
where
\begin{eqnarray*}
q_{p,a,b}(x)&=&x^5-(p+b+a+7)x^4+(bp+ap+4p+ab+5b+5a+18)x^3\\
&-&(abp+2bp+2ap+6p+3ab+8b+8a+22)x^2\\
&+&(bp+ap+4p+2ab+5b+5a+13)x-p-b-a-3.
\end{eqnarray*}
All the eigenvalues larger than $\overline{d}$ are roots of $q_{p,a,b}(x)$. If $x_1 \geq \cdots \geq x_5 > 0$ denote the roots of $q_{p,a,b}(x)$, we know that $S_3(T)=x_1+x_2+x_3=n+4-x_4-x_5$, and hence it suffices to obtain lower bounds for the two smaller roots of $q_{p,a,b}(x)$. We assume that $p \geq 1$, and we shall show later that this implies the general case. We prove that $x_4 \geq \frac{4}{n}$ and $x_5 \geq \frac{2}{n}$, which leads to
$$S_3(T)= n+4-x_4-x_5 < n+4-\frac{6}{n},$$
as desired. To conclude that $x_5 \geq 2/n$, we show that there is a single Laplacian eigenvalue of $T$ that is smaller than $2/n$, which we know to be $0$. To this end, we use the algorithm described in Section~\ref{sec_overview}: given a vertex $v \in V(T)$, the value of $a(v)$ at termination of the algorithm is equal to $a(u)=\frac{n-2}{n}>0$ if $u$ is a leaf. For the roots $v_1$ and $v_2$ of the branches of type 1, we have $a(v_1)=1-\frac{2a}{n-2}-\frac{2}{n}$ and $a(v_2)=1-\frac{2b}{n-2}-\frac{2}{n}$. Note that $a(v_1)$ and $a(v_2)$ cannot both be negative, as otherwise we would have $a,b > \frac{n}{2}-2$, which leads to
$$n=a+b+p+3>\frac{n}{2}-2+\frac{n}{2}-2+p+3=n+p-1\geq n,$$
a contradiction.

We now consider the central vertex $v_0$. Let $$f(p,a,b,n)=a(v_0)=p+2-\frac{2}{n}-\frac{pn}{n-2}-\frac{1}{a(v_1)}-\frac{1}{a(v_2)}.$$ With the change of variables $p=\rho+1$, $b=\beta+1$ and $a=s+\beta+1$, we see that
\begin{equation}\label{form1.2}
n^3(n-2)^3a(v_1)a(v_2)f(\rho+1,s+\beta+1,\beta+1,\rho+2\beta+s+6)
\end{equation}
is a multivariate polynomial on the indeterminates $\rho,\beta,s$ whose coefficients are negative (see equation ~\eqref{form1.2a} in the appendix). Since we know that these indeterminates are nonnegative in our case, this implies that $a(v_1)a(v_2)a(v_0)<0$. Moreover, because $a(v_1)$ and $a(v_2)$ cannot both be negative, we conclude that exactly one of the terms $a(v_1),a(v_2),a(v_0)$ is negative, as required.

On the other hand, to conclude that $x_4 \geq 4/n$, we show that there are at most two Laplacian eigenvalues of $T$ that are smaller than $4/n$. To this end, we look again at the algorithm for eigenvalue location: observe that $a(u)=\frac{n-4}{n}>0$ for every leaf vertex $u$, so there are at most three Laplacian eigenvalues that are smaller than $4/n$. Moreover, this number would be equal to three if and only if $a(v_0),a(v_1),a(v_2)<0$. Assume that $a(v_1),a(v_2)<0$. Since $a(v_1)=1-\frac{4}{n}-\frac{4a}{n-4}$ and $a(v_2)=1-\frac{4}{n}-\frac{4b}{n-4}$, we must have $a,b>\frac{(n-4)^2}{4n}$ and hence
\begin{eqnarray}
&a=n-b-p-3<n-\frac{(n-4)^2}{4n}-p-3, \label{eq1.2aux1}\\
&p=n-a-b-3<n-2\frac{(n-4)^2}{4n}-3=\frac{n}{2}+1-\frac{8}{n}. \label{eq1.2aux2}
\end{eqnarray}
Combining our expression for $a(v_1)$ with~\eqref{eq1.2aux1}, we obtain
\begin{eqnarray*}
a(v_1)&=&1-\frac{4}{n}-\frac{4a}{n-4}>1-\frac{4}{n}-\frac{4n}{n-4}-\frac{n-4}{n}+\frac{4(p+3)}{n-4}\\
&=& -2-\frac{8}{n}+\frac{4(p-1)}{n-4} = -\frac{2(n^2+2n-2np-16)}{n(n-4)},
\end{eqnarray*}
and $a(v_2)$ satisfies the same bound. This will be used to show that $a(v_0)>0$:
\begin{eqnarray*}
a(v_0)&=&p+2-\frac{4}{n}-\frac{pn}{n-4}-\frac{1}{a(v_1)}-\frac{1}{a(v_2)}\\
&>&2-\frac{4}{n}-\frac{4p}{n-4}+\frac{n(n-4)}{n^2+2n-2np-16}.
\end{eqnarray*}
For a fixed $n\geq8$, this expression, viewed as a real function with indeterminate $p$, has two relative extrema, one between 2 and $\frac{n}{2}$ (a relative minimum) and another greater than $\frac{n}{2}+1$ (a relative maximum), namely
\begin{eqnarray*}
j_1(n)=\frac{1}{2}\left(1-\frac{1}{\sqrt{2}}\right)n-\frac{8}{n}+\sqrt{2}+1,\\
j_2(n)=\frac{1}{2}\left(1+\frac{1}{\sqrt{2}}\right)n-\frac{8}{n}-\sqrt{2}+1.
\end{eqnarray*}
Moreover, it is easy to see that $a(v_0)$ achieves an absolute minimum (for a real number $p$) when $p=j_1(n)$, so that
$$a(v_0)\geq \frac{4\,\left( {n}^{2}-4\,\sqrt{2}\,n-4\,n+12\,\sqrt{2}\right) }{\sqrt{2}\,\left( n-4\right) \,n}\geq2\sqrt{2}-\frac{5}{2}>0,$$
as required.

To conclude the proof in this case, observe that for $p=0$, the number $1$ is a root of $q_{0,a,b}(x)$, so that $x_4+x_5 \geq 1+ x_5 >6/n$, as $n \geq 6$ for all trees in this class.

\vspace{5pt}

\noindent \emph{Case 2 $(\sigma=r)$:} Assume that $T$ be an $n$-vertex tree in $\mathfrak{F}_4$ such that $\sigma(T)=r$ and let $\r$ be the number of branches of type 1 for which $s_i \geq 2$. By Lemma~\ref{lemadiam4}, we know that $\r \leq 2$ and that the root of $T$ is not adjacent to pendants. The remainder of the proof depends on the value of $\r$.

\vspace{5pt}

\noindent \emph{Case 2.1 $(\r = 2)$:} Let $a=s_r$ and $b=s_{r-1}$ be the number of leaves of the branches of type 1 with the highest number of leaves. So $n=a+b+2r-1$ and  the characteristic polynomial of the Laplacian matrix of $T$ is given by
\begin{equation}\label{mona6}
P_T(x)=p_{a,b,r}(x)(x^2-3x+1)^{r-3}(x-1)^{a+b-2}x,
\end{equation}
where
\begin{eqnarray*}
p_{a,b,r}(x)&=&x^6-(r+b+a+7)x^5+(br+ar+6r+ab+5b+5a+19)x^4 -\\
&-& (abr+4br+4ar+14r+3ab+9b+9a+24)x^3 +\\
&+& (2abr+5br+5ar+16r+3ab+7b+7a+13)x^2 -\\
&-& (2br+2ar+9r+2ab+3b+3a+1)x+2r+b+a-1
\end{eqnarray*}
The Laplacian eigenvalues of $T$ are $\frac{3\pm\sqrt{5}}{2}$ with multiplicity $r-3$, 0, 1 (with multiplicity $a+b-2$) and the six roots $x_1\geq x_2\geq\ldots\geq x_6>0$ of $p_{a,b,r}(x)$. Since $\sigma=r$, the $r$ largest roots of $P_T(x)$ are $\frac{3+\sqrt{5}}{2}>2$, with multiplicity $r-3$, and the three largest roots of $p_{a,b,r}(x)$. Consider the change of variables $a=\alpha+2$, $b=\beta+2$ and $r=\rho+4$. We may see that
\begin{eqnarray*}
p_{a,b,r}\left(\frac{9}{5}\right)&=&p_{\alpha+2,\beta+2,\rho+4}\left(\frac{9}{5}\right)\\
&=&\frac{5476}{3125}\rho+\frac{5906}{3125}\alpha+\frac{5906}{3125}\beta+\frac{666}{625}\beta \rho+ \frac{666}{625}\alpha \rho+\frac{1071}{625} \alpha \beta+\frac{81}{125} \alpha \beta \rho +\frac{24716}{15625},
\end{eqnarray*}
that is, $p_{a,b,r}\left(\frac{9}{5}\right)$ may be written as a multivariate polynomial with indeterminates $\alpha,\beta,\rho$ with nonnegative coefficients. In particular, $p_{a,b,r}\left(\frac{9}{5}\right)> 0$ for all $a,b \geq 2$ and all $r \geq 4$, and hence $\frac{9}{5}< x_4<2-2/n$. So, for $r\geq4$, we have
    \begin{eqnarray}
    S_{r}(T)&=& (r-3)\frac{3+\sqrt{5}}{2}+x_1+x_2+x_3 \nonumber \\
    &\leq& (r-3)\frac{3+\sqrt{5}}{2}+(a+b+r+7)-x_4  \nonumber\\
    &\leq& \frac{5+\sqrt{5}}{2}r+a+b+\frac{5-3\sqrt{5}}{2}-\frac{9}{5}  \nonumber\\
    &=& \frac{5+\sqrt{5}}{2}r+a+b-\frac{15\sqrt{5}-7}{10} < 4r+a+b-4 \label{wh1} \\
   &=& n+2r-3 < n-2+2r-\frac{2r}{n}. \nonumber
    \end{eqnarray}
The last inequality follows from $n\geq2r+3$, while the validity of~\eqref{wh1} comes from $r\geq4$.

If $r=3$, we proceed similarly. With $a=\alpha+2$ and $b=\beta+2$, we conclude that
$$p_{a,b,3}\left(1+\frac{6}{n}\right)=p_{\alpha+2,\beta+2,3}\left(1+\frac{6}{\alpha+\beta+9}\right)$$
is a rational function on $\alpha$ and $\beta$ whose denominator is positive for $\alpha,\beta \geq 0$ and whose numerator is a bivariate polynomial with nonnegative coefficients. Hence $x_4\geq 1+\frac{6}{n}$, which implies that
\begin{eqnarray*}
S_3(T)&=&x_1+x_2+x_3<a+b+10-x_4\leq a+b+9-\frac{6}{n} \\
&=&n+4-\frac{6}{n} = n-2+2\cdot3-\frac{2\cdot3}{n}.
\end{eqnarray*}
Finally, in the case $r=2$, the polynomial $p_{a,b,2}(x)$ has a factor $(x^2-3x+1)$, which cancels out the missing factor (exponent -1) of
$P_T(x)$ in~\eqref{mona6}. The two largest Laplacian eigenvalues are the two largest roots $x_1 \geq x_2$ of
\begin{eqnarray*}
q_{a,b}(x)=x^4-(b+a+6)x^3+(ab+4b+4a+12)x^2-(2ab+4b+4a+10)x+b+a+3.
\end{eqnarray*}
Mimicking the previous argument, we can see that $x_3$ (the third largest root of $q$) satisfies $1+\frac{4}{n}\leq x_3$ with the change of variables $a=\alpha+2$ and $b=\beta+2$, as
\begin{eqnarray*}
q_{a,b}\left(1+\frac{4}{n}\right)=q_{\alpha+2,\beta+2}\left(1+\frac{4}{\alpha+\beta+7}\right)
\end{eqnarray*}
is a rational function on $\alpha$ and $\beta$ whose denominator is positive for $\alpha,\beta \geq 0$ and whose numerator has negative coefficients. This implies that
\begin{eqnarray*}
S_2(T)&=&x_1+x_2=a+b+6-x_3-x_4<n+3-1-\frac{4}{n} \\
&=&n-2+2\cdot2-\frac{2\cdot2}{n},
\end{eqnarray*}
as required.

\vspace{5pt}

\noindent \emph{Case 2.2 $(\r = 1)$:}  Consider a tree $T \in \mathfrak{F}_4$ with $r$ branches and $s_r=s \geq 2$, while $s_i=1$ for $i \in \{1,\ldots,r-1\}$. The characteristic polynomial of the Laplacian matrix of $T$ may be derived from the previous case by setting $a=s$ and $b=1$, which leads to
$$P_T(x)=p_{s,r}(x)(x^2-3x+1)^{r-2}(x-1)^{s-1}x,$$ where
$$p_{s,r}(x)=x^4-(s+r+5)x^3+(rs+3s+4r+8)x^2-(2rs+2s+5r+4)x+s+2r.$$
The largest eigenvalues of $T$ are $\frac{3+\sqrt{5}}{2}$, with multiplicity $r-2$, and the largest roots $x_1 \geq x_2$ of $p_{s,r}(x)$. As $p_{s,r}(1)=-(r-1)s$, the smallest root $x_4$ of $p_{s,r}(x)$ is inside $(0,1)$. Also, since $p_{s,r}(2)=s$, we know that $x_3\in(1,2)$. We look for better lower bounds on these two roots. If we consider the parameters $s=\alpha+2$ and $r=\beta+3$ (or $r=2$), we obtain
\begin{eqnarray*}
p_{s,2}\left(\frac{3}{2(s+4)}\right)&=&p_{\alpha+2,2}\left(\frac{3}{2(\alpha+6)}\right)=\\
&=&\frac{16\alpha^5+336\alpha^4+2724\alpha^3+10818\alpha^2+22086\alpha+20493}{16(\alpha+6)^4} > 0    \\
p_{s,2}\left(1+\frac{5}{2(s+4)}\right)&=&p_{\alpha+2,2}\left(1+\frac{5}{2(\alpha+6)}\right)=\\
&=&-\frac{16\alpha^5+376\alpha^4+3184\alpha^3+11606\alpha^2+16334\alpha+4427}{16(\alpha+6)^4} < 0
\end{eqnarray*}
\begin{equation*}
p_{s,r}\left(\frac{33}{20}\right)=p_{\alpha+2,\beta+3}\left(\frac{33}{20}\right)=-\frac{92400\alpha\beta+161140\beta+57140\alpha+22279}{20^4} < 0
\end{equation*}
For $r=2$, we have $x_4>\frac{3}{2n}$ and $x_3>1+\frac{5}{2n}$, while,  for $r\geq3$, we have $x_3>\frac{33}{20}$.

If $r=2$, the number of vertices is $n=s+4$. We obtain
\begin{eqnarray*}
S_2(T)&=&x_1+x_2=s+7-(x_3+x_4) \\
&<& s+7-\left(1+\frac{4}{n}\right)=n+2-\frac{4}{n}.
\end{eqnarray*}
If $r \geq 3$, we have $n=2r+s$, so that
\begin{eqnarray*}
S_r(T)&=&x_1+x_2+(r-2)\frac{3+\sqrt{5}}{2}=s+r+5-x_3-x_4+(r-2)\frac{3+\sqrt{5}}{2} \\
&<&\frac{5+\sqrt{5}}{2}r+s+2-\frac{33}{20}-\sqrt{5} \\
&<&4r+s-3=n+2r-3<n-2+2r-\frac{2r}{n}.
\end{eqnarray*}

\vspace{5pt}

\noindent \emph{Case 2.3 $(\r = 0)$:} The characteristic polynomial of the Laplacian matrix is equal to $$P_T(x)=(x^2-(r+3)x+2r+1)(x^2-3x+1)^{r-1}x.$$ The $r$ eigenvalues larger than $\overline{d}$ are $\frac{3+\sqrt{5}}{2}$ with multiplicity $r-1$ and $\frac{r+3+\sqrt{r^2-2r+5}}{2}$, which is bounded above by $r+1+\frac{1}{r-1}$. We have
\begin{eqnarray*}
S_r(T) &<& r+1+\frac{1}{r-1}+(r-1)\frac{3+\sqrt{5}}{2} \\
&=& \frac{5+\sqrt{5}}{2}r-\frac{1+\sqrt{5}}{2}+\frac{1}{r-1} \\
&\leq& \frac{5+\sqrt{5}}{2}r-\frac{1+\sqrt{5}}{2}+\frac{1}{2} \\
&=& \frac{5+\sqrt{5}}{2}r-\frac{\sqrt{5}}{2} < 4r-2
\end{eqnarray*}
for $r>\frac{4-\sqrt{5}}{3-\sqrt{5}}\approx 2.309$. We may suppose that $r \geq 3$, as we would have $T=\mathcal{P}_5$ for $r=2$. Using $n=2r+1$ we get
\begin{eqnarray*}
S_r(T) &<& 4r-2=2r-1+2r-1 \\
&=& n-2+2r-1 < n-2+2r-\frac{2r}{2r+1} \\
&=&n-2+2r-\frac{2r}{n},
\end{eqnarray*}
as required. This concludes the proof of Theorem~\ref{teo1} for $\mathfrak{F}_4$.
\end{proof}

\subsection{Trees in $\mathfrak{F}_5$}
In this section we will consider the family $\mathfrak{F}_5$, whose elements are $n$-vertex trees with diameter five consisting of a root vertex adjacent to $p \geq 0$ pendant vertices (branches of type 0), to $r-1 \geq 1$ branches of type 1 (where the $i$th branch has $s_i$ leaves) and to exactly one branch of type 2, which has $t$ leaves. In particular, we have $n=p+r+2+t+\sum_{i=1}^{r-1} s_i$.
As in the case of trees in $\mathfrak{F}_4$, an important ingredient in proving Theorem~\ref{teo1} is to estimate the number of Laplacian eigenvalues of an $n$-vertex tree that are larger than the average $\overline{d}=2-2/n$.
\begin{Lema}{\label{lemadiam5}} \cite{FHRT11} A tree in $\mathfrak{F}_5$ with $r-1 \geq 1$ branches of type 1 has $r+1$ or $r+2$ eigenvalues larger than the average $\overline{d}$.
    \end{Lema}

\begin{proof}[Proof of Theorem~\ref{teo1} for $\mathfrak{F}_5$] There are two basic cases in this proof. If $r \geq 4$, we use the decomposition of Lemma~\ref{aplic_decomp_legal}, while, for $r \in \{2,3\}$, we analyze the characteristic polynomial directly.

\vspace{5pt}
\noindent \emph{Case 1 $(r \geq 4)$:} For a tree $T$ in $\mathfrak{F}_5$ with $r\geq4$, we use
the decomposition of Lemma~\ref{aplic_decomp_legal} applied to $T$ with $k \in \{r+1,r+2\}$ to obtain $$S_{k}(T) \leq S_{k+1}(A),$$
where $A$ is the matrix defined in (\ref{eqndiamgeralcp}) or (\ref{eqndiamgeralsp}) according to whether $v_0$ is adjacent to pendant vertices (i.e. $A=A_p$) or not (i.e. $A=A_{\bar{p}}$). The $r+2$ largest eigenvalues of $A_p$ and $A_{\bar{p}}$ are the same in this case, namely the degree $\delta=r+p$ of the root, the largest eigenvalue of each block $T_i$ and the two largest eigenvalues of the single block $Q$ associated with a brach of type 2 (both are larger than 1 by Lemma~\ref{lemamatriza3}). Moreover, if we combine Lemma~\ref{decomp_Rojo} with Lemmas~\ref{lemamatriza2} and \ref{lemamatriza3}, we see that the $(r+3)$-rd largest eigenvalue of $A$ is at most 1 (it is exactly 1 if there are pendants).

Let $y_3$ denote the smallest eigenvalue of $Q$. In the case $k=r+1$,  we see that
    \begin{eqnarray}
    S_{r+1}(T)&<& \delta+\sum_{i=1}^{r-1}\left(s_i+2-\frac{1}{s_1+2}\right)+t+4-y_3 \nonumber\\
    &=&r+p+\sum_{i=1}^{r-1} s_i +t+4+2(r-1)-\sum_{i=1}^{r-1}\frac{1}{s_i+2}-y_3 \nonumber\\
    &=&n+2r-\sum_{i=1}^{r-1}\frac{1}{s_i+2}-y_3.
    \end{eqnarray}
To conclude the proof, we show that
\begin{equation}\label{mona4}
\sum_{i=1}^{r-1}\frac{1}{s_i+2}+y_3 \geq \frac{2r+2}{n}.
\end{equation}
To this end, we note that, for $S=\max \{s_i \colon i=1 \ldots,r-1\}$, we have
$$\sum_{i=1}^{r-1}\frac{1}{s_i+2} \geq \max\left\{ \frac{r-1}{S+2},\frac{(r-1)^2}{n+r-t-4}\right\}.$$
The first lower bound comes from our choice of $S$, and, for the second lower bound, we used Lemma \ref{lemacotasomainvsi2} with $\sum s_i = n-r-t-p-2$.

If $t=1$, we use the lower bound $y_3>0.19 \geq \frac{1}{n}$ (as $n \geq 10$) to conclude that
\begin{eqnarray*}
&\frac{(r-1)^2}{n+r-t-4} + y_3 \geq \frac{2r+2}{n}\\
 \Longrightarrow & \frac{(r-1)^2}{n+r-t-4} \geq \frac{2r+1}{n}\\
 \Longleftrightarrow & (n-2)r^2 - (4n-2t-7)r+t+4 \geq 0.
\end{eqnarray*}
The two roots of $p(r)= (n-2)r^2 - (4n-2t-7)r+t+4$ are nonnegative and smaller than $\frac{4n-2t-7}{n-2}<4$, which implies our result. For $t \geq 2$, the same conclusion would be reached if $t \leq \frac{n}{4}$, since we have $y_3 > \frac{1}{4t}  \geq \frac{1}{n}$ in this case, and the previous argument would hold. So, we assume that $t>n/4$.

We shall prove \eqref{mona4} by showing that
$$ \frac{r-1}{S+2} \geq \frac{2r+2}{n}.$$
Note that the condition $r \geq 4$ implies that
$$S \leq \frac{1}{3} \left(n-r-t-2\right) < \frac{1}{3} \left(n-r-\frac{n}{4}-2\right)=\frac{n}{4}-\frac{r+2}{3}. $$
In particular,
$$ \frac{r-1}{S+2}  >  \frac{r-1}{n/4-(r+2)/3+2} \geq \frac{4(r-1)}{n} \geq \frac{2r+2}{n},$$
as required. Therefore, the upper bound of Theorem~\ref{teo1} is true for $k=r+1$.

For $k=r+2$, we recall that the $(r+3)$-rd largest eigenvalue of $A$ is at most 1, so that, applying the above bound, we get
    \begin{eqnarray*}
    S_{r+2}(T)&\leq&S_{r+1}(T)+1\\
    &<&n+2r-\frac{2(r+1)}{n} + 1\\
    &=&n+2(r+1)-\frac{2(r+2)}{n}-\left(1-\frac{2}{n} \right)<n-2+2(r+2)-\frac{2(r+2)}{n},
    \end{eqnarray*}
given that $n > 2$.

\vspace{5pt}

\noindent \emph{Case 2 $(r=2)$:} In this case, we look at the characteristic polynomial of the Laplacian matrix directly. Assuming that the branches of type 1 and 2 have $s$ and $t$ leaves, respectively, we have $P_T(x)=q_{s,t,p}(x)\cdot(x-1)^{s+t+p-3}\cdot x,$ where
    \begin{eqnarray*}
    q_{s,t,p}(x)&=&x^6-(t+s+p+9)x^5+(st+pt+7t+ps+7s+6p+31)x^4\\
    &-&(pst+5st+4pt+17t+4ps+17s+13p+53)x^3\\
    &+&(2pst+7st+5pt+18t+4ps+18s+13p+48)x^2\\
    &-&(3st+2pt+8t+ps+8s+6p+22)x+t+s+p+4.
    \end{eqnarray*}
Lemma~\ref{lemadiam5} implies that $3\leq\sigma\leq4$. Using the eigenvalue localization algorithm mentioned in Section~\ref{sec_overview}, we may see that there are three eigenvalues larger than 2 and one eigenvalue inside $(1,2)$. Moreover, at most three eigenvalues are smaller than 1 (one of them is 0). All eigenvalues that are larger than $\overline{d}$ are roots of $q_{s,t,p}(x)$, which we denote $x_1 \geq x_2 \geq \cdots \geq x_6$. We claim that $x_6 \geq 2/n$, $x_5 \geq 4/n$ and $x_4 \geq 7/5$, and the strategy to prove this resembles what has been done when considering trees in $\mathfrak{F}_4$ (see Case 1.2).

As before, we first assume that $p \geq 1$. To show that $x_6=1/n$, we apply the localization algorithm on $T$ with $\alpha=\frac{1}{n}$. We obtain $a(u)=1-\frac{1}{n}>0$ for all leaf vertices $u$. As $n=p+s+t+4$, the maximum value for $s$ and $t$ is $n-5$. For the remaining four vertices, the algorithm gives us
\begin{eqnarray*}
&& a(v)=s+1-\frac{1}{n}-\frac{ns}{n-1}\geq\frac{4}{n-1}-\frac{1}{n}>0,\\
&& a(w)=t+1-\frac{1}{n}-\frac{nt}{n-1}>0,\\
&& a(w')=2-\frac{1}{n}-\frac{1}{a(w)},\\
&& a(v_0)=p+2-\frac{1}{n}-\frac{pn}{n-1}-\frac{1}{a(v)}-\frac{1}{a(w')}.
\end{eqnarray*}
To analyze the sign of $a(w')$ and $a(v_0)$, we calculate
\begin{equation}\label{whoopie1}
n^3(n-1)^3a(v)a(w)a(v')a(v_0),
\end{equation}
where we replace $n$ by $p+s+t+4$. This gives a multivariate polynomial with indeterminates $p,s,t$ and negative coefficients (see equation~\eqref{whoopie1a} in the appendix). In particular, exactly one of $a(w')$ and $a(v_0)$ is negative, and hence there is only one eigenvalue smaller than $\frac{1}{n}$, which we know to be 0. This implies that $x_6 \geq 1/n$.

To see that $x_5 \geq 2/n$, we apply the algorithm with $\alpha=\frac{2}{n}$ to conclude that $a(u)$ is nonnegative for all leaf vertices and that $a(v)$ and $a(w)$ cannot both be negative. (The calculations are very similar to the ones used for establishing $x_4 \geq 4/n$ when considering Case 1.2 for trees in $\mathfrak{F}_4$.) For $w'$ and $v_0$, we obtain
\begin{eqnarray*}
&& a(w')=2-\frac{2}{n}-\frac{1}{a(w)}\\
&& a(v_0)=p+2-\frac{2}{n}-\frac{pn}{n-2}-\frac{1}{a(v)}-\frac{1}{a(w')}.
\end{eqnarray*}
We claim that at most two of the values $a(v_0),a(v),a(w),a(w')$ may be negative, which implies the desired result. If $a(v)$ and $a(w)$ are both positive, we are done, so suppose that this is not the case. If $a(w)<0$ (and hence $a(v) \geq 0$), then $a(w')$ must be positive and we are done. Finally, suppose that $a(v)<0$ (and hence $a(w) \geq 0$). If $a(v') \geq 0$, we are done, while $a(v')<0$ leads to $$a(v_0)>p+2-\frac{2}{n}-\frac{pn}{n-2}=2-\frac{2}{n}-\frac{2p}{n-2}\geq\frac{4}{n-2}-\frac{2}{n}>0$$ because $p\leq n-6$, implying our claim.

Finally, to prove that $x_4\geq 7/5$, we apply the algorithm for $\alpha=\frac{7}{5}$. We obtain $a(u)=-\frac{2}{5}$ whenever $u$ is a leaf vertex. For $v$ and $w$, we obtain $a(v)=s+1-\frac{7}{5}+\frac{5s}{2}=\frac{35s-4}{10}>0$ and $a(w)=\frac{35t-4}{10}>0$, while, for $w'$ we have $a(w')=2-\frac{7}{5}-\frac{10}{35t-4}=\frac{105t-62}{5(35t-4)}>0$. Moreover, we have
\begin{eqnarray*}
a(v_0)&=&p+2-\frac{7}{5}+\frac{5p}{2}-\frac{1}{a(v)}-\frac{1}{a(w')}\\
&=&\frac{7\,\left( 18375\,\rho\,\tau\,\varsigma+12775\,\tau\,\varsigma+7525\,\rho\,\varsigma+1065\,\varsigma+16275\,\rho\,\tau+9815\,\tau+6665\,\rho+329\right) }{10\,\left( 105\,\tau+43\right) \,\left( 35\,\varsigma+31\right) }
\end{eqnarray*}
where $\rho=p+1$, $\varsigma=s+1$ and $\tau=t+1$. As a consequence, we have four eigenvalues greater than $\frac{7}{5}$, which implies $x_4>\frac{7}{5}$.

Recall that $\sigma \in \{3,4\}$. Using the relation $\sum_{i=1}^6 x_i=t+s+p+9=n+5$, we have
\begin{eqnarray*}
S_3(T) &=& x_1+x_2+x_3=n+5-x_{4}-x_{5}-x_6\\
&<& n+5-\frac{7}{5}-\frac{2}{n}-\frac{1}{n}<n+4-\frac{6}{n}
\end{eqnarray*}
for $n\geq8$. In the case $n=7$, we verify directly that $S_3(T)<10<11-\frac{6}{7}$ for the single tree $T$ in $\mathfrak{F}_5$ with at least one pendant vertex at the root. For $\sigma=4$ we have
\begin{eqnarray*}
S_4 &=& x_1+x_2+x_3+x_4=n+5-x_{5}-x_6\\
&<& n+5-\frac{3}{n}<n+6-\frac{8}{n}.
\end{eqnarray*}

To conclude the proof, we need to consider the case $p=0$, where the result holds easily since we may easily show that $x_4=x_5=1$ (for instance, by showing that 1 is a root of multiplicity two of $q_{s,t,0}(x)$).

\vspace{5pt}

\noindent \emph{Case 3 $(r=3)$:} We proceed as in the case $r=2$. Indeed, we may show that the characteristic polynomial of the Laplacian matrix has a single root smaller than $1/n$ (which we know to be 0) and at most three roots smaller than $2/n$. We may also show that $\mu_5 \geq \frac{7}{5}$. As before, the first conclusion may be obtained by applying the localization algorithm for $\alpha=\frac{1}{n}$: we conclude that $a(v_1),a(v_2),a(w)>0$ and that $a(u)>0$ for every leaf vertex. Moreover, the product
\begin{equation*}
n^5(n-1)^4a(v_1)a(v_2)a(w)a(w')a(v_0)
\end{equation*}
is always negative, so that exactly one eigenvalue is smaller than $1/n$. The second conclusion may be obtained easily, since, by the localization algorithm, at most one among $a(v_1),a(v_2),a(w)$ may be negative. To reach the third conclusion, we note that the application of the localization algorithm leads to $a(u)=-\frac{2}{5}$ for all leaf vertex $u$, while $a(v_1)=\frac{7s_1}{2}-\frac{2}{5}>0$, $a(v_2)=\frac{7s_2}{2}-\frac{2}{5}>0$ and $a(w)=\frac{7t}{2}-\frac{2}{5}>0$. We also obtain $a(w')=\frac{105t-62}{5(35t-4)}>0$. Finally, we note that
\begin{eqnarray*}
a(v_0) &=&7\,(643125\,\rho\,\tau\,\varsigma_1\,\varsigma_2+630875\,\tau\,\varsigma_1\,\varsigma_2+263375\,\rho\,\varsigma_1\,\varsigma_2+112525\,\varsigma_1\,\varsigma_2+569625\,\rho\,\tau\,\varsigma_2\\
&+&506275\,\tau\,\varsigma_2+233275\,\rho\,\varsigma_2+78165\,\varsigma_2+569625\,\rho\,\tau\,\varsigma_1+506275\,\tau\,\varsigma_1+233275\,\rho\,\varsigma_1+78165\,\varsigma_1\\
&+&504525\,\rho\,\tau+401915\,\tau+206615\,\rho+50189)/(10\,\left( 105\,\tau+43\right) \,\left( 35\,\varsigma\_1+31\right) \,\left( 35\,\varsigma\_2+31\right) )
\end{eqnarray*}
where $p=\rho+1$, $s_i=\varsigma_i+1$ and $t=\tau+1$. This is positive, so that $|\{v \colon a(v) > 0\}| \geq 5$ and therefore $\mu_5 > \frac{7}{5}$.

The Laplacian spectrum of $T$ contains one eigenvalue zero, $n-9$ eigenvalues 1, three eigenvalues in $(0,1]$, one in $(1,2)$ and four greater than 2. As $\sum_{i=1}^n \mu_i = 2n-2$, we have
\begin{eqnarray*}
S_4(T) &=& \mu_1+\mu_2+\mu_3+\mu_4=2n-2-(n-9)\cdot1-\mu_{5}-\mu_{n-3}-\mu_{n-2}-\mu_{n-1}\\
&<& n+7-\frac{7}{5}-\frac{2}{n}-\frac{1}{n}-\frac{1}{n}\leq n+6-\frac{8}{n}
\end{eqnarray*}
for $n\geq10$. For $n=9$, we verify this inequality directly for every such tree (there is a single tree with this property. which satisfies $S_4(T) <14 < 15-8/9$). For $\sigma=5$ we have
\begin{eqnarray*}
S_5 &=& \mu_1+\mu_2+\mu_3+\mu_4+\mu_5\\
&=&2n-2-(n-9)\cdot1-\mu_{n-3}-\mu_{n-2}-\mu_{n-1}\\
&<& n+7-\frac{4}{n}<n+8-\frac{10}{n}.
\end{eqnarray*}
Finally, if $p=0$, we have
\begin{eqnarray*}
S_4 &=& \mu_1+\mu_2+\mu_3+\mu_4=2n-2-(n-9)\cdot1-\mu_{5}-\mu_{n-3}-\mu_{n-2}-\mu_{n-1}\\
&<& n+7-1-1-\frac{2}{n}<n+6-\frac{8}{n}
\end{eqnarray*}
for every $n\geq8$, as we may show that $\mu_5=\mu_{n-3}=1$.
\end{proof}

\subsection{Trees in $\mathfrak{F}_6$}
To conclude this section, we consider the family $\mathfrak{F}_6$ of trees with diameter 6, which consist of a root vertex adjacent to $p \geq 0$ pendants, $r_1 \geq 0$ branches of type 1 and $r_2 \geq 2$ branches of type 2. As with the special families with diameter four and five, there is information about the number of Laplacian eigenvalues that are larger than the average.
\begin{Lema}\cite{FHRT11}{\label{lemadiam6}} Given $T \in \mathfrak{F}_6$ with $r$ branches of type 1 and 2, the number of Laplacian eigenvalues larger than the average $\overline{d}$ is between $r+1$ and $r+4$.
\end{Lema}

We are ready to prove Theorem~\ref{teo1} for trees in $\mathfrak{F}_6$.
\begin{proof}[Proof of Theorem~\ref{teo1} for $\mathfrak{F}_6$]  Let $T$ be a tree in $\mathfrak{F}_6$ with $r_1$ branches of type 1 and $r_2$ branches of type 2 and, as usual, let $s_i$ and $t_j$ be the number of leaves on the $i$-th branch of type 1 and on the $j$-th branch of type 2, respectively, ordered in nondecreasing order. Once again, the idea is to evaluate the sum of eigenvalues by decomposing the matrix associated with $T$ to obtain $$S_{k}(T) \leq S_{k+1}(A),$$
where $A$ is the matrix defined in (\ref{eqndiamgeralcp}) or (\ref{eqndiamgeralsp}) according to whether $v_0$ is adjacent to pendant vertices (i.e. $A=A_p$) or not (i.e. $A=A_{\bar{p}}$). By Lemma~\ref{lemadiam6}, it suffices to consider $k=r + \ell \in \{r+1,\ldots,r+4\}$.

By Lemmas~\ref{lemamatriza2} and~\ref{lemamatriza3}, the $r+1$ largest eigenvalues of $A$ are
$\delta$ and the largest eigenvalue associated with each $T_i$ and $Q_j$. Following these, we have the second largest eigenvalues of the $Q_j$, the eigenvalue 1 (when there are pendants), and finally the smallest eigenvalues of $T_i$ and $Q_j$ (which are smaller than 1). For the $r_1$ largest eigenvalues in matrices associated with branches of type 1, we shall use the upper bound $s_i+2-\frac{1}{s_i+2}$ (the $i$-th branch of this type has $s_i$ leaves). Analogously, for the $r_2$ largest eigenvalues in matrices associated with branches of type 2, we shall use the upper bound $t_j+2+\frac{1}{4t_j}$ (the $j$-th branch of this type has $t_j$ leaves).

Recall that $n=1+p+r_1+2r_2+\sum_{r_1}s_i+\sum_{r_2}t_j$. By Lemma~\ref{lemamatriza3} (the sum of the eigenvalues of $Q_j$ is $t_j+4$), the addition of the second, or even the third largest eigenvalue associated with a branch of type 2 would increase the upper bound of the previous paragraph by at most $2-\frac{1}{4t_j} \geq \frac{7}{4}$, as $t_j \geq 1$. On the other hand, the addition of an eigenvalue 1 or a smallest eigenvalue of a matrix associated with a branch of type 1 contributes with at most 1 in this upper bound (the contribution of an additional eigenvalue of the branch with type 1 is at most $\frac{1}{s_i+2}$). In other words, for $\ell \leq r_2$, we have
\begin{eqnarray}
S_{r+\ell}(T) & \leq & \delta + \sum_{i=1}^{r_1} \left(s_i+2-\frac{1}{s_i+2}\right) + \sum_{j=1}^{r_2} \left(t_j+2+\frac{1}{4t_j}\right)+ \sum_{j=1}^{\ell} \left(2-\frac{1}{4t_{r_2-j+1}}\right) \nonumber\\
&=& \delta+\sum_{i=1}^{r_1}s_i+\sum_{j=1}^{r_2}t_j + 2r_1+2r_2+ 2 \ell+ \sum_{j=1}^{r_2-\ell}\frac{1}{4t_j}-\sum_{i=1}^{r_1}\frac{1}{s_i+2}. \label{mona2}
\end{eqnarray}

We shall consider the cases $k=r+1$ and $k \geq r+2$ separately. First assume that $k=r+1$. We have (recall that $s_1\leq s_i$)
    \begin{eqnarray*}
    S_{r+1}(T) &<& \delta+\sum_{i=1}^{r_1}s_i + \sum_{i=1}^{r_2}t_j + 2r_1+2r_2+2-\sum_{i=1}^{r_1}\frac{1}{s_i+2} + \sum_{j=1}^{r_2-1}\frac{1}{4t_j}\\
    &\leq& (n-1-r_2)+2r+2-\frac{r_1}{s_1+2}-\sum_{j=1}^{r_2-1}\left(1-\frac{1}{4t_j}\right)+r_2-1 \\
    &\leq& n-2+2(r+1)-\frac{r_1}{s_1+2}-(r_2-1)\frac{3}{4}
    \end{eqnarray*}
To complete the proof, we show that $$\frac{r_1}{s_1+2}+\frac{3r_2-3}{4}\geq\frac{2r+2}{n}.$$
If $r_1=0$, it is easy to see that, since $r_2 \geq 2$, the above inequality is satisfied when $n\geq8\geq\frac{8}{3}\frac{r_2+1}{r_2-1}$. There is a single tree in $\mathfrak{F}_6$ with fewer than eight vertices, namely the path $\mathcal{P}_7$, for which the result may be verified directly. In fact, $\sigma=3$ in this case and $S_3(\mathcal{P}_7)<10,$ which is smaller than the upper bound given in~\eqref{new_bound}.

If $r_1\geq1$, we have $n\geq3r_2+(s_1+1)r_1+1$. With some algebraic manipulations, we confirm that $$(3r_2+(s_1+1)r_1+1)(4r_1+(3r_2-3)(s_1+2))\geq8(r_1+r_2+1)(s_1+2)$$
for $r_1\geq1$, $r_2\geq2$ and $s_1\geq1$. Hence, our result follows from
    \begin{eqnarray*}
    n(4r_1+(3r_2-3)(s_1+2))&\geq&8(r_1+r_2+1)(s_1+2)\\
    \frac{4r_1+(3r_2-3)(s_1+2)}{4(s_1+2)}&\geq&\frac{2(r_1+r_2+1)}{n}\\
    \frac{r_1}{s_1+2}+\frac{3r_2-3}{4}&\geq&\frac{2r+2}{n}.
    \end{eqnarray*}

We move to the case $k=r+\ell>r+1$. From equation~\eqref{mona2}, we derive the following for $\ell \leq r_2$:
\begin{eqnarray*}
 S_{r+\ell}(T) &<& n-1+ (2(r+\ell) - 1)-r_2+1+\frac{r_2-\ell}{4}-\sum_{i=1}^{r_1}\frac{1}{s_i+2}\\
&=& n-1+ (2(r+\ell) - 1)-\left(\frac{3r_2}{4}+\frac{\ell}{4}-1-\sum_{i=1}^{r_1}\frac{1}{s_i+2}\right).
\end{eqnarray*}
Moreover, if $r_2 < \ell$, say $r_2=\ell-\ell^{\ast}$, an upper bound is obtained when the contribution of $\ell^{\ast}$ of the additional eigenvalues is set to 1, from which one may easily derive
$$S_{r+\ell}(T) \leq  n-1+ (2(r+\ell) - 1)-\left(2\ell-\ell^{\ast}-1\right).$$
However, note that
$$2\ell-\ell^{\ast} \geq \frac{3r_2}{4}+\frac{\ell}{4} = \ell-\frac{3 \ell^{\ast}}{4},$$
since this is equivalent to the inequality $\ell \geq \frac{\ell^{\ast}}{4}$, which holds trivially.
Hence, to conclude our proof, it suffices to show that
\begin{equation}\label{eqfrida}
\frac{3r_2}{4}+\frac{\ell}{4}-1+\sum_{i=1}^{r_1}\frac{1}{s_i+2} > \frac{2(r+\ell)}{n},
\end{equation}
which is equivalent to
\begin{equation*}
n > \frac{8(r+\ell)}{3r_2+\ell-4+\sum_{i=1}^{r_1}\frac{1}{s_i+2}}
=8\left(1+\frac{r_1-2r_2+4-\sum_{i=1}^{r_1}\frac{1}{s_i+2}}{3r_2+\ell-4+\sum_{i=1}^{r_1}\frac{1}{s_i+2}} \right).
\end{equation*}
Since $r_2 \geq 2$, we have
$$\frac{r_1-2r_2+4-\sum_{i=1}^{r_1}\frac{1}{s_i+2}}{3r_2+\ell-4+\sum_{i=1}^{r_1}\frac{1}{s_i+2}} \leq \frac{r_1-\sum_{i=1}^{r_1}\frac{1}{s_i+2}}{4+\sum_{i=1}^{r_1}\frac{1}{s_i+2}}.$$
In particular,
\begin{equation}\label{mona3}
8\left(1+\frac{r_1-2r_2+4-\sum_{i=1}^{r_1}\frac{1}{s_i+2}}{3r_2+\ell-4+\sum_{i=1}^{r_1}\frac{1}{s_i+2}} \right)\leq 8 + \frac{8r_1-\sum_{i=1}^{r_1}\frac{1}{s_i+2}}{4+\sum_{i=1}^{r_1}\frac{1}{s_i+2}}\leq 8+2r_1.
\end{equation}
Note that the last inequality may be replaced by a strict inequality unless $r_1=0$.
Since $T$ contains $r_1$ branches of type 1 (each of which has at least two vertices) and $r_2 \geq 2$ branches of type 2 (each of which contains at least three vertices), we know that $n \geq p + 1 + 2r_1 + 3r_2 \geq 7+2r_1$, which, in light of~\eqref{mona3}, verifies equation~(\ref{eqfrida}) for all cases other than $T=\mathcal{P}_7$, for which $r_2=2,r_1=0,p=0,t_1=t_2=1.$ However, for $T=\mathcal{P}_7$, we know that $\sigma=3=r+1$, and hence this case does not apply.
\end{proof}

\section{Proof of Theorem~\ref{teo1} - General Case}\label{sec_general}

In this section, we give a full proof of Theorem \ref{teo1} using the results of the previous sections.

\begin{Lema}{\label{lemamu1naoestrela}}\cite{FHRT11} If $T\not = \mathcal{S}_n$ is a tree with $n$ vertices and $\mu_1$ is the largest Laplacian eigenvalue of $T$, then  $$\mu_1<n-\frac{1}{2}.$$
\end{Lema}

As a reminder, we restate our main result.\\
\textbf{Theorem \ref{teo1}:} \emph{Every tree $T$ with $n \geq 6$ vertices and
diameter greater than or equal to four
satisfies}
\begin{equation}\label{tapsi2}
S_k(T) < n-2+2k-\frac{2k}{n}.
\end{equation}

Observe that the single tree with $n \leq 5$ vertices and diameter at least four is the path $\mathcal{P}_5$. Unfortunately, this tree does not satisfy~\eqref{tapsi2} for $k=2$, as
$$S_2(\mathcal{P}_5)=4+2\left(\cos(\pi/5)+\cos(2\pi/5) \right) >\frac{31}{5}.$$
On the other hand, it is easy to check that the bound holds for every $k \neq 2$, and that we may use the upper bound $S_2(\mathcal{P}_5) < \frac{25}{4}.$

\begin{proof}[Proof of Theorem~\ref{teo1}] For a contradiction, we suppose that the theorem fails, and we fix an $n$-vertex counterexample $T$ with the smallest number of vertices. In light of our work in Section~\ref{sec_particular}, we know that $T \notin \mathfrak{F}_4 \cup \mathfrak{F}_5 \cup \mathfrak{F}_6$. In particular, we have $n \geq 7$, as there are precisely three nonisomorphic trees with diameter larger than or equal to four, one of which belongs to $\mathfrak{F}_5$, and two of which belong to $\mathfrak{F}_4$.

The remainder of our argument is based on choosing an edge $e \in E$ whose endpoints have degree larger than one, so as to consider the two components $T_1$ and $T_2$ of the forest $F=T-e$, both of which contain an edge. The existence of such an edge is guaranteed by the condition on the diameter of $T$. Moreover, we may suppose that $e$ was chosen with the additional property that the two components are not stars, because, if we suppose that the removal of every such edge $e$ of $T$ produces a star, then $T \in \mathfrak{F}_4 \cup \mathfrak{F}_5 \cup \mathfrak{F}_6$, which we know not to be the case.

Let $n_1$ and $n_2$ denote the number of vertices in $T_1$ and $T_2$, respectively, so that $n_1+n_2=n$. \\
For a convenient labelling of the vertices, the Laplacian matrices of $T$ and $F$ satisfy
    \begin{equation}{\label{teoeqmtmf}}
    M_T=M_F+M=M_F+\begin{pmatrix}
    M_{*} & 0 \\ 0 & 0
    \end{pmatrix},
    \end{equation}
where $M_{*}=\begin{pmatrix} 1 & -1 \\ -1 & 1 \end{pmatrix}$. The spectrum of $M$ has one element equal to 2, while all the others are 0. By Theorem \ref{teoautmat} (with index set $I=\{1,\ldots,k\}$) we obtain the relation
    \begin{eqnarray}\label{tapsi1}
    \sum_{i=1}^{k}\lambda_i(M_T) &\leq& \sum_{i=1}^{k}\lambda_i(M_F) + \sum_{i=1}^{k}\lambda_i(M_{*}) \nonumber \\
    S_k(T) &\leq& S_k(F)+2.
    \end{eqnarray}
Since $F$ is disconnected, its eigenvalues are precisely the eigenvalues of $T_1$ and $T_2$. Of the $k$ largest Laplacian eigenvalues of $F$, we assume that $k_1$ are eigenvalues of $T_1$ and $k_2$ are eigenvalues of $T_2$, with $k_1+k_2=k$. The inequality~\eqref{tapsi1} may be rewritten as
    \begin{equation}{\label{eqnsk1t1k2t2}}
    S_k(T)\leq S_{k_1}(T_1)+S_{k_2}(T_2)+2.
    \end{equation}
Without loss of generality, we suppose that $T_1$ has diameter larger than or equal to the diameter of $T_2$. We shall now split our argument into a few cases, according to the diameters of $T_1$ and $T_2$:
    \begin{enumerate}
    \item $T_1$ and $T_2$ have diameter 3;
    \item $T_1$ and $T_2$ have diameter $\geq$ 4;
    \item  $T_1$ has diameter $\geq 4$ and $T_2$ has diameter 3;
    \end{enumerate}

\vspace{5pt}

\noindent \textbf{Case 1:} In the case when $T_1$ and $T_2$ have diameter three, we apply the upper bound of Lemma~\ref{teodiam3}(c) to the inequality (\ref{eqnsk1t1k2t2}). Each $T_i$ has exactly two eigenvalues larger than the average, so that $\sigma(F) = 4$. Because the Laplacian eigenvalues do not decrease with the addition of a new edge (and hence $\mu_i(T)\geq\mu_i(F)$ for every $i$), we conclude that $\sigma(T)\geq4$. In particular, Lemma~\ref{lemassigma} tells us that (\ref{eqnsk1t1k2t2}) must fail for some $k \geq 4$.

However, for $k=4$, we have
    \begin{eqnarray*}
    S_4(T)&<& n_1+2-\frac{2}{n_1}+n_2+2-\frac{2}{n_2}+2\\
    &=&n+6-\frac{2}{n_1}-\frac{2}{n_2}=n+6-\frac{2n}{n_1n_2}\\
    &\leq&n+6-\frac{8}{n}= (n-2)+2 \cdot 4 - \frac{2 \cdot 4}{n},
    \end{eqnarray*}
so that~(\ref{tapsi2}) holds for $T$ in this case. For $k\geq5$ (remembering that $k_1,k_2\geq2$), we get
    \begin{eqnarray*}
    S_k(T)&<& n_1+k_1-\frac{2}{n_1}+n_2+k_2-\frac{2}{n_2}+2\\
    &=&n+k+2-\frac{2}{n_1}-\frac{2}{n_2}<n-2+2k-(k-4)\\
    &\leq&n-2+2k-\frac{2k}{n},
    \end{eqnarray*}
because $n\geq6\geq\frac{2k-4}{k-4}$, so that~(\ref{tapsi2}) holds for $T$ in this case. As a consequence, a minimum counterexample to the statement of the theorem cannot qualify for Case 1.

\vspace{5 pt}

\noindent \textbf{Case 2:}
We consider the case when $T_1$ and $T_2$ have diameter larger than or equal to four. Due to the minimality of $T$, we know that both $T_1$ and $T_2$ satisfy~(\ref{tapsi2}), unless $T_i=\mathcal{P}_5$ and $k_i=2$ for some $i \in \{1,2\}$.

First suppose that this is not the case. We have
    \begin{eqnarray*}
    S_k(T)&<& n_1-2+2k_1-\frac{2k_1}{n_1} + n_2-2+2k_2-\frac{2k_2}{n_2} +2\\
    &=&n-2+2k-\frac{2k_1}{n_1}-\frac{2k_2}{n_2} \leq n-2+2k-\frac{2k_1+2k_2}{n_1+n_2}\\
    &=&n-2+2k-\frac{2k}{n},
    \end{eqnarray*}
as the inequality $\frac{a+b}{c+d} \leq \frac{a}{c}+\frac{b}{d}$ holds for every $a,b \in \mathbb{N}$ and $c,d \in \mathbb{N}_{>0}$. Thus (\ref{tapsi2}) holds for $T$ in this case.

We now assume that $T_2=\mathcal{P}_5$ and $k_2=2$, while this is not the case for $T_1$. For $\mathcal{P}_5$, we use the upper bound $S_2(\mathcal{P}_5)<\frac{25}{4}$. Observe that $k=k_1+2$ and $n=n_1+5$. It follows that
    \begin{eqnarray*}
    S_k(T)&=& S_{k_1+2} < n_1-2+2k_1-\frac{2k_1}{n_1}+S_2(\mathcal{P}_5)+2\\
    &<&n_1+2k_1-\frac{2k_1}{n_1}+\frac{25}{4}=(n_1+5)-2+2(k_1+2)-\frac{2k_1}{n_1}-\frac{3}{4}\\
    &=&n-2+2k-\frac{2k_1}{n_1}-\frac{3}{4}<n-2+2k-\frac{2k_1}{n}-\frac{4}{n}=n-2+2k-\frac{2k}{n},
    \end{eqnarray*}
because $n\geq 10$, which again tells us that $T$ does not contradict~(\ref{tapsi2}).

Finally, if $T_1=T_2=\mathcal{P}_5$ and $k_1=k_2=2$, we have
    \begin{eqnarray*}
    S_4(T)&\leq& S_2(\mathcal{P}_5)+S_2(\mathcal{P}_5)+2<\frac{25}{4}+\frac{25}{4}+2=\frac{29}{2}\\
    &<&\frac{76}{5}=10-2+8-\frac{8}{10}.
    \end{eqnarray*}
This shows that $T$ cannot qualify for Case 2.

\vspace{5pt}

\noindent \textbf{Case 3:} We assume that $T_1$ has diameter larger than or equal to four and $T_2$ has diameter three. In the following, we shall use~(\ref{tapsi2}) for $T_1$, unless $T_1=\mathcal{P}_5$ and $k_1=2$, in which case we use $S_2(\mathcal{P}_5)<\frac{25}{4}$. The upper bound on $S_{k_2}(T_2)$ depends on the value of $k_2$.

If $k_2=0$ (recall that $n\geq n_1+4 \geq 9$), we have
    \begin{eqnarray*}
    S_k(T)&\leq& S_k(T_1)+2<n_1-2+2k-\frac{2k}{n_1}+2\\
    &=&(n_1+4)-2+2k-\frac{2k}{n_1}<n-2+2k-\frac{2k}{n},\\
    S_2(T)&\leq& S_2(\mathcal{P}_5)+2<\frac{25}{4}+2=\frac{33}{4}\\
    &<&\frac{95}{9}\leq n+2-\frac{4}{n}.
    \end{eqnarray*}
If $k_2=1$, we use the upper bound on $S_1(T_2)$ of Lemma~\ref{lemamu1naoestrela}, so that
    \begin{eqnarray*}
    S_k(T)&\leq& S_{k-1}(T_1)+S_1(T_2)+2 < n_1-2+2(k-1)-\frac{2k-2}{n_1}+n_2-\frac{1}{2}+2\\
    &=& n-2+2k-\frac{2k-2}{n_1}-\frac{1}{2}<n-2+2k-\frac{2k-2}{n}-\frac{2}{n}= n-2+2k-\frac{2k}{n},
    \end{eqnarray*}
    \begin{eqnarray*}
    S_3(T)&\leq& S_2(\mathcal{P}_5)+S_1(T_2)+2<\frac{25}{4}+n_2-\frac{1}{2}+2\\
    &=& (n_2+5)+4-\frac{5}{4}\leq n+4-\frac{6}{n}.
    \end{eqnarray*}
Finally, if $k_2\geq2$, we bound $S_{k_2}(T_2)$ by means of Lemma~\ref{teodiam3}(c), obtaining, for $T_1 \neq \mathcal{P}_5$,
    \begin{eqnarray*}
    S_k(T)&<& n_1-2+2k_1-\frac{2k_1}{n_1}+n_2+k_2-\frac{2}{n_2}+2\\
    &=& n-2+2k-(k_2-2)-\frac{2k_1}{n_1}-\frac{2}{n_2}.
    \end{eqnarray*}
The validity of~\eqref{tapsi2} is a consequence of $\frac{2k_1}{n_1}+\frac{2}{n_2}>\frac{2k}{n}$, which follows from
$$\frac{2k_1}{n_1}+\frac{2}{n_2}-\frac{2k}{n}=\frac{1}{n_1n_2n}\left((k_2-2)n_1n_2(n-2)+(2k_1-1)n_2^2+n_1^2+(n_1-n_2)^2\right)>0$$
If $T_1 = \mathcal{P}_5$, we get
    \begin{eqnarray*}
    S_k(T)&\leq& S_2(\mathcal{P}_5)+S_{k-2}(T_2)+2<\frac{25}{4}+n_2+k-2-\frac{2}{n_2}+2\\
    &=& (n_2+5)-2+2k+\frac{13}{4}-\frac{2}{n_2}-k < n-2+2k-\frac{2k}{n}
    \end{eqnarray*}
because $n_2<n$ and $n\geq9>\frac{8(k-1)}{4k-13}$ for $k\geq4$. In particular, no counterexample to~\eqref{tapsi2} can qualify for Case 3. This implies that there is no such counterexample, concluding the proof of Theorem~\ref{teo1}.
\end{proof}

\bibliographystyle{acm}
\bibliography{ref}

\appendix

\section{Additional proofs}

In this appendix, we include the proofs that have been omitted in the main text. 

\subsection{Proof of Theorem~\ref{order}} To conclude the proof of Theorem~\ref{order} for $n$ even, we need to show that $\left\lfloor \sqrt{n-2-\frac{8}{n}+\frac{16}{n^2}} \right\rfloor =  \left\lfloor \sqrt{n-3} \right\rfloor$. Clearly, we have $\lfloor \sqrt{n-3} \rfloor \leq \left\lfloor \sqrt{n-2-\frac{8}{n}+\frac{16}{n^2}} \right\rfloor\leq \lfloor \sqrt{n-2} \rfloor$, since $n-3\leq n-2-\frac{8}{n}+\frac{16}{n^2}\leq n-2$ for $n\geq4$. Fix nonnegative integers $q$ and $m$ such that $n=q^2+m$ and $q$ is largest with this property (in particular, $m \leq 2q$). It is easy to see that $\lfloor \sqrt{n-3} \rfloor = \lfloor \sqrt{n-2} \rfloor = q-1$ if $m \in \{0,1\}$ and  $\lfloor \sqrt{n-3} \rfloor = \lfloor \sqrt{n-2} \rfloor = q$ if $m \in \{3,\ldots,2q\}$. We now suppose that $m=2$, in which case $$q-1=\lfloor \sqrt{q^2-1} \rfloor\leq \left\lfloor \sqrt{q^2-\frac{8}{q^2+2}+\frac{16}{(q^2+2)^2}}\right\rfloor.$$
Let $h=\frac{8}{q^2+2}-\frac{16}{(q^2+2)^2}$, which satisfies $0<h<1<q^2$ if $n \geq 6$. The Taylor series expansion of $\sqrt{q^2-h}$ yields
\begin{eqnarray*}
    \sqrt{q^2-h} &=& \sqrt{q^2}+\frac{1(-h)}{1!2\sqrt{q^2}}-\frac{1(-h)^2}{2!4\sqrt{q^2}^3}+\frac{3(-h)^3}{3!8\sqrt{q^2}^5}-\frac{3\cdot5(-h)^4}{4!16\sqrt{q^2}^7}+\cdots \\
    &=&q-\frac{h}{2q}-\sum_{i=2}^{\infty}\frac{(2i-3)!h^i}{(i-2)!i!2^{2i-2}q^{2i-1}} <q-\frac{h}{2q}.
\end{eqnarray*}
Therefore $\lfloor\sqrt{q^2-h}\rfloor\leq\lfloor q -\frac{h}{2q} \rfloor<q$, which implies that $\left\lfloor\sqrt{n-2-\frac{8}{n}+\frac{16}{n^2}}\right\rfloor=q-1$. Similar calculations lead to desired result for odd values of $n$.

\subsection{Proof of Theorem~\ref{thm_counter}}  According to Lemma~\ref{teodiam3}, the characteristic polynomial of the Laplacian matrix of $T(n-3,1)$ is
$$P(x)=(x^3-nx^2-2x^2+3nx-2x-n)\cdot(x-1)^{n-4}\cdot x.$$
Moreover, there are exactly two roots larger than $\overline{d}=2-2/n$, which are the largest roots of the factor $p_{n-3,1}(x)$ of degree three. We may see, as a consequence of $p_{n-3,1}(\ell_n)<0$,  that the smallest root of this factor, which we call $x_3$, satisfies $x_3 > \ell_n=\frac{3-\sqrt{5}}{2}+\frac{0.4}{n}$.

The formula for the characteristic polynomial of the Laplacian matrix of $F(n,k)$ depends on the remainder of $n$ modulo three:
    $$q(x)=q_{n~(\mbox{mod } 3)}(x)\cdot(x^2-(k+1)x+1)\cdot(x-1)^{n-6}\cdot x,$$ where
    \begin{eqnarray*}
    && q_0(x)=(x^2-(k+3)x+3)\cdot(x-k) \\
    && q_1(x)=x^3-(2k+4)x^2+(k+2)^2x-3k-1 \\
    && q_2(x)=x^3-(2k+5)x^2+(k^2+5k+5)x-3k-2.
    \end{eqnarray*}
Let us compute the difference between the Laplacian energies of $F(n,k)$ and $T(n-3,1)$. Consider $\Delta_{LE} = LE(F(n,k)) - LE(T(n-3,1))$. We know that the number of Laplacian eigenvalues larger than the average degree $\overline{d}=2-2/n$ is equal to three for $F(n,k)$ (see Lemma~\ref{lemadiam4}) and equal to two for $T(n-3,1)$ (see Lemma~\ref{teodiam3}). For $n=3k$, we have, using Proposition~\ref{prop_1},
    \begin{eqnarray*}
    \Delta_{LE} &=& 2\left[k+\frac{k+3+\sqrt{k^2+6k-3}}{2}+\frac{k+1+\sqrt{k^2+2k-3}}{2}\right]-6\overline{d} -2(n+2-x_3)+4\overline{d}\\
    &>&2\left[k+\left(k+3-\frac{5}{2k}-\frac{1}{100}\right)+\left(k+1-\frac{1}{k}\right)\right]-2\left[3k+2-\frac{3-\sqrt{5}}{2}-\frac{0.4}{3k}\right]-2\overline{d}\\
    &=&4-\frac{1}{50}-\frac{7}{k}+3-\sqrt{5}+\frac{0.8}{3k}-2\left(2-\frac{2}{3k}\right)\\
    &=&-\frac{5.4}{k}+3-\sqrt{5}-\frac{1}{50}>0
    \end{eqnarray*}
for $k  \geq 8$. For $k=6,7$, it can be seen in Table \ref{tableasmall}.

For $n=3k+1$, considering $y_3$ the smallest root of $q_1(x)$, we have
\begin{eqnarray*}
    \Delta_{LE} &=& 2\left[2k+4-y_3+\frac{k+1+\sqrt{k^2+2k-3}}{2}\right]-6\overline{d}-2(n+2-x_3)+4\overline{d}\\
    &>& 2\left[2k+4-y_3+\left(k+1-\frac{1}{k}\right)\right]-2\left[3k+1+2-\frac{3-\sqrt{5}}{2}-\frac{0.4}{3k+1}\right]-2\left(2-\frac{2}{3k+1}\right)\\
    &=& -2y_3-\frac{2}{k}+3-\sqrt{5}+\frac{4.8}{3k+1}
\end{eqnarray*}
For $\Delta_{LE}>0$, it is enough to have $y_3<\frac{3-\sqrt{5}}{2}+\frac{2.4}{3k+1}-\frac{1}{k}$. In fact, when we evaluate $q_1$ at this value, we obtain a positive number for $k\geq6$ (this value is greater than the smaller root), and since this has been checked for $k=5$ (see Table \ref{tableasmall}), the result follows.

For $n=3k+2$, we have
\begin{eqnarray*}
    \Delta_{LE} &=& 2\left[2k+5-y_3+\frac{k+1+\sqrt{k^2+2k-3}}{2}\right]-6\overline{d}-2(n+2-x_3)+4\overline{d}\\
    &>& 2\left[2k+5-y_3+\left(k+1-\frac{1}{k}\right)\right]-2\left[3k+2+2-\frac{3-\sqrt{5}}{2}-\frac{0.4}{3k+2}\right]-2\left(2-\frac{2}{3k+2}\right)\\
    &=& -2y_3-\frac{2}{k}+3-\sqrt{5}+\frac{4.8}{3k+2}
\end{eqnarray*}
To have $\Delta_{LE}>0$, it suffices to verify that $y_3<\frac{3-\sqrt{5}}{2}+\frac{2.4}{3k+2}-\frac{1}{k}$. We evaluate $q_2$ at this value and obtain a positive number for $k\geq6$. As the case $k=5$ has been checked directly (see Table~\ref{tableasmall}), the result follows.

\subsection{Additional expressions}
The polynomial defined in equation~\eqref{form1.2} is such that
\begin{equation}\label{form1.2a}
\begin{split}
\frac{-1}{2}& n^3(n-2)^3a(v_1)a(v_2) f(\rho+1,s+\beta+1,\beta+1,\rho+2\beta+s+6)=2080+1792 s+4672\rho\\
&+3584\beta+3336\rho^2+8256\rho\beta+4128\rho s+2304\beta^2+2304\beta s+1080 s^2+2664\rho^2\beta s\\
&+3312 \rho\beta^2 s+2024\rho\beta s^2+5984\rho\beta s+1256\rho^3\beta+704\rho^2 s^2+2392\rho^2 s+2664\rho^2 \beta^2\\
&+4784 \rho^2\beta+628\rho^3 s+460\rho s^3+1736\rho s^2+2208\beta^3\rho+5984\beta^2\rho+128\beta^3 s+680\beta^2 s^2\\
&+1392\beta s^2+8\rho^5\beta+28\rho^4\beta^2+160\rho^4\beta +80\rho^4 s+116\rho^3 s^2+108\rho^2 s^3+56\rho^3\beta^3\\
&+456\rho^3\beta^2+4\rho^5 s+7\rho^4 s^2+234\rho^4+1192\rho^3+150 s^4+536 s^3+\rho^6+24\rho^5+s^6\\
&+20 s^5+28\rho^4\beta s+84\rho^3\beta^2 s+44\rho^3\beta s^2+456\rho^3\beta s+108\rho^2\beta^2s^2+44\rho^2\beta s^3\\
&+552\rho^2\beta s^2+76\rho\beta^2 s^3+28\rho\beta s^4+376\rho\beta s^3+128\rho^2\beta^3 s+80\beta^4\rho s+104\beta^3\rho s^2\\
&+792\beta^2\rho s^2+1008\rho^2\beta^2 s+832\beta^3\rho s+8\rho^3 s^3+7\rho^2 s^4+4\rho s^5+68\rho s^4+64\beta^4\rho^2\\
&+672\beta^3\rho^2+16\beta^4 s^2+32\beta^3 s^3+160\beta^3 s^2+24\beta^2s^4+240\beta^2s^3+ 8\beta s^5+120\beta s^4\\
&+64\beta^4+640\beta^3+32\beta^5\rho+416\beta^4\rho+960\beta^2 s+616\beta s^3.
\end{split}
\end{equation}

The polynomial defined in equation~\eqref{whoopie1} is given by
\begin{equation}\label{whoopie1a}
\begin{split}
&n^3(n-1)^3a(v)a(w)a(v')a(v_0)=\\
&-({t}^{7}+4\,s\,{t}^{6}+5\,p\,{t}^{6}+20\,{t}^{6}+6\,{s}^{2}\,{t}^{5}+18\,p\,s\,{t}^{5}+67\,s\,{t}^{5}+11\,{p}^{2}\,{t}^{5}+87\,p\,{t}^{5}+172\,{t}^{5} \nonumber\\
&+5\,{s}^{3}\,{t}^{4}+28\,p\,{s}^{2}\,{t}^{4}+88\,{s}^{2}\,{t}^{4}+38\,{p}^{2}\,s\,{t}^{4}+273\,p\,s\,{t}^{4}+487\,s\,{t}^{4}+15\,{p}^{3}\,{t}^{4}+170\,{p}^{2}\,{t}^{4}\nonumber\\
&+651\,p\,{t}^{4}+839\,{t}^{4}+5\,{s}^{4}\,{t}^{3}+32\,p\,{s}^{3}\,{t}^{3}+82\,{s}^{3}\,{t}^{3}+64\,{p}^{2}\,{s}^{2}\,{t}^{3}+403\,p\,{s}^{2}\,{t}^{3}+601\,{s}^{2}\,{t}^{3}\nonumber\\
&+52\,{p}^{3}\,s\,{t}^{3}+531\,{p}^{2}\,s\,{t}^{3}+1809\,p\,s\,{t}^{3}+2049\,s\,{t}^{3}+15\,{p}^{4}\,{t}^{3}+210\,{p}^{3}\,{t}^{3}+1128\,{p}^{2}\,{t}^{3}+2753\,p\,{t}^{3}\nonumber\\
&+2566\,{t}^{3}+6\,{s}^{5}\,{t}^{2}+33\,p\,{s}^{4}\,{t}^{2}+88\,{s}^{4}\,{t}^{2}+74\,{p}^{2}\,{s}^{3}\,{t}^{2}+441\,p\,{s}^{3}\,{t}^{2}+601\,{s}^{3}\,{t}^{2}+84\,{p}^{3}\,{s}^{2}\,{t}^{2}\nonumber\\
&+798\,{p}^{2}\,{s}^{2}\,{t}^{2}+2460\,p\,{s}^{2}\,{t}^{2}+2420\,{s}^{2}\,{t}^{2}+48\,{p}^{4}\,s\,{t}^{2}+625\,{p}^{3}\,s\,{t}^{2}+3057\,{p}^{2}\,s\,{t}^{2}+6653\,p\,s\,{t}^{2}\nonumber\\
&+5418\,s\,{t}^{2}+11\,{p}^{5}\,{t}^{2}+180\,{p}^{4}\,{t}^{2}+1198\,{p}^{3}\,{t}^{2}+4072\,{p}^{2}\,{t}^{2}+7081\,p\,{t}^{2}+5034\,{t}^{2}+4\,{s}^{6}\,t\nonumber\\
&+22\,p\,{s}^{5}\,t+67\,{s}^{5}\,t+53\,{p}^{2}\,{s}^{4}\,t+330\,p\,{s}^{4}\,t+487\,{s}^{4}\,t+72\,{p}^{3}\,{s}^{3}\,t+683\,{p}^{2}\,{s}^{3}\,t+2097\,p\,{s}^{3}\,t\nonumber\\
&+2049\,{s}^{3}\,t+58\,{p}^{4}\,{s}^{2}\,t+739\,{p}^{3}\,{s}^{2}\,t+3489\,{p}^{2}\,{s}^{2}\,t+7197\,p\,{s}^{2}\,t+5418\,{s}^{2}\,t+26\,{p}^{5}\,s\,t\nonumber\\
&+414\,{p}^{4}\,s\,t+2635\,{p}^{3}\,s\,t+8389\,{p}^{2}\,s\,t+13354\,p\,s\,t+8484\,s\,t+5\,{p}^{6}\,t+95\,{p}^{5}\,t+756\,{p}^{4}\,t\nonumber\\
&+3241\,{p}^{3}\,t+7928\,{p}^{2}\,t+10524\,p\,t+5931\,t+{s}^{7}+6\,p\,{s}^{6}+20\,{s}^{6}+16\,{p}^{2}\,{s}^{5}+106\,p\,{s}^{5}+172\,{s}^{5}\nonumber\\
&+25\,{p}^{3}\,{s}^{4}+246\,{p}^{2}\,{s}^{4}+795\,p\,{s}^{4}+839\,{s}^{4}+25\,{p}^{4}\,{s}^{3}+324\,{p}^{3}\,{s}^{3}+1560\,{p}^{2}\,{s}^{3}+3297\,p\,{s}^{3}\nonumber\\
&+2566\,{s}^{3}+16\,{p}^{5}\,{s}^{2}+256\,{p}^{4}\,{s}^{2}+1630\,{p}^{3}\,{s}^{2}+5160\,{p}^{2}\,{s}^{2}+8105\,p\,{s}^{2}+5034\,{s}^{2}+6\,{p}^{6}\,s\nonumber\\
&+114\,{p}^{5}\,s+900\,{p}^{4}\,s+3785\,{p}^{3}\,s+8952\,{p}^{2}\,s+11292\,p\,s+5931\,s+{p}^{7}+22\,{p}^{6}+207\,{p}^{5}\nonumber\\
&+1083\,{p}^{4}+3413\,{p}^{3}+6498\,{p}^{2}+6939\,p+3213)
\end{split}
\end{equation}

\end{document}